\documentclass[11pt]{article}
\usepackage{amssymb,amsmath,amsthm, amsfonts}

\def\sqr#1#2{{\vcenter{\vbox{\hrule height.#2pt
              \hbox{\vrule width.#2pt height#1pt \kern#1pt \vrule width.#2pt}
          \hrule height.#2pt}}}}
\def\signed #1{{\unskip\nobreak\hfil\penalty50
          \hskip2em\hbox{}\nobreak\hfil#1
          \parfillskip=0pt \finalhyphendemerits=0 \par}}
\def\endpf{\signed {$\sqr69$}}
\def\sqr#1#2{{\vcenter{\vbox{\hrule height.#2pt
              \hbox{\vrule width.#2pt height#1pt \kern#1pt \vrule width.#2pt}
              \hrule height.#2pt}}}}
\def\signed #1{{\unskip\nobreak\hfil\penalty50
              \hskip2em\hbox{}\nobreak\hfil#1
              \parfillskip=0pt \finalhyphendemerits=0 \par}}
\def\endpf{\signed {$\sqr69$}}
\def\3n{\negthinspace \negthinspace \negthinspace }
\def\2n{\negthinspace \negthinspace }
\def\1n{\negthinspace }

\def\={\buildrel \triangle \over =}

\def\ms{\medskip}

\def\esssup{\mathop{\rm esssup}}

\def\max{\mathop{\rm max}}
\def\min{\mathop{\rm min}}
\def\exp{\mathop{\rm exp}}
\def\sup{\mathop{\rm sup}}

\def\inf{\hbox{\rm inf$\,$}}

\def\|{\Big |}
\def\({\Big (}
\def\){\Big )}
\def\[{\Big[}
\def\]{\Big]}
\def\be{\begin{equation}}
\def\bel{\begin{equation}\label}
\def\ee{\end{equation}}
\def\bt{\begin{theorem}}
\def\bcd{\begin{condition}}
\def\ecd{\end{condition}}
\def\et{\end{theorem}}
\def\bc{\begin{corollary}}
\def\ec{\end{corollary}}
\def\bde{\begin{definition}}
\def\ede{\end{definition}}
\def\bl{\begin{lemma}}
\def\el{\end{lemma}}
\def\bp{\begin{proposition}}
\def\ep{\end{proposition}}
\def\br{\begin{remark}}
\def\er{\end{remark}}
\def\ba{\begin{array}}
\def\ea{\end{array}}
\def\ed{\end{document}}

\def\square#1{\vbox{\hrule\hbox{\vrule height#1%
     \kern#1\vrule}\hrule}}
\def\rectangle#1#2{\vbox{\hrule\hbox{\vrule height#1%
     \kern#2\vrule}\hrule}}

\font\tenbb=msbm10 \font\sevenbb=msbm7 \font\fivebb=msbm5

\newfam\bbfam
\scriptscriptfont\bbfam=\fivebb \textfont\bbfam=\tenbb
\scriptfont\bbfam=\sevenbb

\newtheorem{lemma}{Lemma}[section]
\newtheorem{remark}{Remark}[section]

\newtheorem{theorem}{Theorem}[section]
\newtheorem{corollary}{Corollary}[section]

\newtheorem{definition}{Definition}[section]
\newtheorem{proposition}{Proposition}[section]
\newtheorem{condition}{Condition}[section]

\makeatletter
   
   \@addtoreset{equation}{section}
\makeatother

\begin{document}
\title{Optimal Stochastic Control with Recursive Cost Functionals
 of Stochastic Differential Systems Reflected in a Domain}

\author{Juan Li\thanks{School of Mathematics and Statistics, Shandong University (Weihai), Weihai 264200, P. R. China.
This author has been supported by the NSF of P.R.China (Nos. 11071144, 11171187, 11222110), Shandong Province (Nos. BS2011SF010, JQ201202), SRF for ROCS (SEM), Program for New Century Excellent Talents in University (No. NCET-12-0331), 111 Project (No. B12023). {\small\it E-mail:} {\small\tt juanli@sdu.edu.cn}.\ms} \quad and \quad Shanjian Tang\thanks{ Institute of Mathematics
and Department of Finance and Control Sciences, School of Mathematical Sciences, Fudan University, Shanghai 200433, China, and Graduate
Department of Financial Engineering, Ajou University, San 5, Woncheon-dong, Yeongtong-gu, Suwon, 443-749, Korea. This author is supported in
part by NSFC Grant \# 11171076, by Basic Research Program of China (973 Program) Grant \# 2007CB814904, by the Science Foundation of the
Ministry of Education of China Grant \#200900071110001, and by WCU (World Class University) Program through the Korea Science and Engineering
Foundation funded by the Ministry of Education, Science and Technology (R31-20007).  {\small\it E-mail:} {\small\tt sjtang@fudan.edu.cn}.\ms}}

\date{August 22, 2013}

\maketitle

\abstract{ In this paper we study the optimal stochastic control
problem for stochastic differential systems reflected in a domain.
The cost functional is a recursive one, which is defined via
generalized backward stochastic differential equations developed by
Pardoux and Zhang~\cite{PZ}. The value function is shown to be the unique viscosity solution to the associated Hamilton-Jacobi-Bellman equation, which is a fully nonlinear parabolic partial differential equation with a nonlinear Neumann boundary condition. For this, we also prove some new estimates for stochastic differential systems reflected in a domain.}



\newpage
\section{Introduction}
\ \ \ \ \ Let D be an open connected
bounded convex subset of ${\mathbb{R}}^d$ such that
 $D=\{\phi>0\},\ \partial D=\{\phi=0\}$ for some function $\phi\in
 C_b^2({\mathbb{R}}^d)$ satisfying  $|\nabla\phi(x)|=1$ at any $x\in \partial D.$ Note that at any $x\in \partial D$, $\nabla\phi(x)$ is a unit normal vector on the boundary point $x$, pointing towards the interior of $D$.

 Let U  be a metric space.
An admissible control process  is a $U$-valued
${\mathbb{F}}$-progressively measurable process. The set of all
admissible control processes is denoted by ${\mathcal{U}}$. In this paper, for the initial data $(t, x)\in [0, T]\times {\mathbb{R}}^d $ we consider the optimal control problem for the following stochastic differential equations (SDEs) reflected on domain $D$:
 \be\label{0.1}
\left\{
  \begin{array}{rcl}
  X_s&=& x+\int_t^sb(r,X_r,u_r)dr+\int_t^s\sigma(r,X_r, u_r)\, dB_r
 +\int_t^s\nabla\phi(X_r)\, dK_r,\ s\in[t,T];\\
  K_s&=& \int_t^sI_{\{X_r\in\partial D\}}dK_r,\ \ K \mbox{ \rm is
  increasing}.
  \end{array}
  \right.
  \ee
Here, $u(\cdot)\in {\mathcal{U}}$ is an admissible control,  and  the drift $b:[0,T]\times {\mathbb{R}}^d\times U \rightarrow {\mathbb{R}}^d$
  and the diffusion $\sigma: [0,T]\times {\mathbb{R}}^d\times U \rightarrow {\mathbb{R}}^{d\times d}$ are assumed to be uniformly Lipschitz continuous and to have a linear growth in the state variable $x$. In view of Proposition 5.1 in the appendix, the above reflected SDE~(\ref{0.1}) has a unique  strong solution for any $u(\cdot)\in {\mathcal{U}}$, which will be denoted by  $(X^{t,x;u}, K^{t,x;u})$.

Then we consider the following  controlled generalized backward stochastic differential equation (GBSDE) where $(X^{t,x;u}, K^{t,x;u})$ is the solution of above reflected SDE (1.1) :
\be\label{0.2}
   \left \{\begin{array}{rcl}
   -dY_s & = & \displaystyle f(s,X^{t,x;u}_s, Y_s, Z_s,
                              u_s)\, ds+ g(s,X^{t,x;u}_s, Y_s)\, dK_s^{t,x;u}-Z_s\, dB_s,\quad s\in [0,T);\\
      Y_T  & = &\displaystyle  \Phi (X^{t,x;u}_T).
   \end{array}\right.
   \ee
Under suitable conditions on the functions $f,\ g$\ and $\Phi$ (see (H3.2) in Section 3 for more details), it has a unique adapted solution (see Pardoux and Zhang~\cite{PZ}), denoted by  $(Y^{t,x;u}, Z^{t,x;u})$ hereafter.  Our optimal control problem is to maximize the cost functional $J(t,x;u):=Y_t^{t,x;u}$ over all admissible controls $u\in \mathcal{U}$. The associated Hamilton-Jacobi-Bellman (HJB) equation turns out to have a nonlinear Neumann boundary condition, and reads as follows:
\be
 \left \{\begin{array}{ll}
 &\!\!\!\!\! \frac{\partial }{\partial t} W(t,x) +  H(t, x, W, DW,
 D^2W)=0, \quad\quad  \hfill  (t,x)\in [0,T)\times {D} ,  \\
 &\!\!\!\!\! \frac{\partial }{\partial n} W(t,x)+g(t,x,W(t,x))=0, \hfill 0\leq t< T, x\in \partial{D};\\
 &\!\!\!\!\! W(T,x) =\Phi (x), \hfill   x \in \bar{D},
 \end{array}\right.
\ee  where at a point $x\in \partial D$, $\frac{\partial }{\partial
n}=\sum_{i=1}^{d}\frac{\partial }{\partial x_i}\phi(x)\frac{\partial
}{\partial x_i}$, and the Hamiltonian $H$ is given by  $$ H(t, x, y, p, A):= \sup_{u
\in U}\{\frac{1}{2}{\text{tr}}(\sigma\sigma^{T}(t, x,
 u)A)+ \langle p, b(t, x, u)\rangle + f(t, x, y, p.\sigma,
u)\},$$  $\mbox{where}\ (t,x,y,p,A)\in [0, T]\times {\mathbb{R}}^n\times {\mathbb{R}}\times {\mathbb{R}}^d\times
{\mathbf{S}}^d$.  We aim at showing that the value function of our optimal control problem is the unique viscosity solution to above HJB equation (1.3).

The linear BSDEs was studied by Bismut in 1973 (see Bismut~\cite{Bismut1,Bismut2,Bismut3}), and the general nonlinear version was studied by Pardoux and
Peng~\cite{PaPe} in 1990. Since then BSDE has received an extensive attention both in the theory and in the application.  The reader is referred to, among others, El Karoui, Peng and Quenez~\cite{ElPeQu}, Darling and Pardoux~\cite{DP}, Pardoux and Peng~\cite{PP}, Peng~\cite{Pe1,Pe2}, Hu~\cite{H}, and Delbaen and
Tang~\cite{DT}. Stochastic differential equations reflected in a domain are referred to Lions~\cite{L}, Lions and Sznitman~\cite{LS}, Menaldi~\cite{M},
Pardoux and Williams~\cite{PW}, Saisho~\cite{S}, among others. Pardoux and Zhang~\cite{PZ} studied BSDEs~\eqref{0.2}, and gave a
probabilistic formula for the solution of a system of parabolic or elliptic semi-linear partial differential equation (PDE) with a nonlinear
Neumann boundary condition. There are also many other works on a PDE with a nonlinear Neumann boundary condition, for example, Boufoussia and Van Casterenb~\cite{BC} gave an approximation result to semilinear parabolic PDEs with Neumann boundary conditions with the help of BSDEs; Day~\cite{Da} studied the  Neumann boundary conditions for viscosity solutions of Hamilton-Jacobi equations. Different from those works we want to study the optimal control problem for stochastic differential systems reflected in a domain, to give the stochastic representation for the solution of HJB equation (1.3) with a nonlinear Neumann boundary condition.

In this paper, the generalized BSDE formulation of dynamic programming  given by Peng~\cite{Pe1,Pe2} for optimally controlled SDEs, is extended to our controlled stochastic differential systems reflected in a domain. The arguments of Buckdahn and Li~\cite{BL} is also generalized to show that our value function $W$ (see (\ref{3.9})) is deterministic (see Proposition 3.1). Since now the associated BSDE is also driven by the increasing process which incorporates the reflection of the system state on the boundary of the given domain, we have many new difficulties, for example, we have to prove the increasing process $K$ satisfies a new important estimate (Proposition 5.3), and  also prove that, under standard assumptions the value of the system path $Y$ at the initial time has linear growth and is locally Lipschitz in the initial (random) state (Proposition 5.2) which improves the estimates on GBSDE of Pardoux and Zhang~\cite{PZ}. Then we can prove the continuity of the value function (Theorem 3.2), and the value function is the unique viscosity solution of the associated HJB equation subject to a nonlinear Neumann boundary condition (Theorem 4.1). On the other hand, with the help of Proposition 3.1 it allows us to prove the dynamic programming principle (DPP in short, see Theorem 3.1) in a straight forward way by adapting to GBSDEs the method of stochastic backward semigroups introduced by Peng~\cite{Pe1}. Furthermore, our proof of Theorem 4.1 differs heavily from the counterpart of either Buckdahn and Li~\cite{BL} or Peng~\cite{Pe1}, the proof becomes more technical due to the Neumann boundary condition. For more details, the reader is referred to among others Lemmas 4.2 and 4.3 and the constructions of BSDEs (\ref{4.10}), (\ref{4.12}), (\ref{4.21}) and (\ref{4.22}), etc. In particular, unlike~\cite{BL} or~\cite{Pe1}, in our context the coefficients are not necessarily continuous in the control variable $u$, and the control $u$ may take values in a possibly noncompact space $U$.

The rest of our paper is organized as follows. In Section 2, we recall some preliminary theory of BSDEs and GBSDEs. In Section 3, we formulate  our optimal stochastic control problem and define the value function $W$. We prove that $W$ is deterministic and satisfies the DPP. Furthermore, we prove that $W$ is continuous. In Section 4, we prove that $W$ is the unique viscosity solution to the associated HJB equation with a nonlinear Neumann boundary condition. In the end, we give some basic properties on GBSDEs associated with forward reflected SDEs in the Appendix (Section 5.1), where Propositions 5.2 and 5.3 contain new results on GBSDEs. For the reader's convenience, the proofs of Proposition 3.1 and Theorem 3.1 are given in Section 5.2.

\section{Preliminaries}

\ \ \ \ We consider the Wiener space $(\Omega, {\cal{F}}, P)$, where $\Omega$
is the set of continuous functions from [0, T] to ${\mathbb{R}}^d$
starting from 0 ($\Omega= C_0([0, T];{\mathbb{R}}^d)$), $ {\cal{F}}
$ the completed Borel $\sigma$-algebra over $\Omega$, and P the
Wiener measure. Let B be the canonical process:
$B_s(\omega)=\omega_s,\ s\in [0, T],\ \omega\in \Omega$. By
${\mathbb{F}}=\{{\mathcal{F}}_s,\ 0\leq s \leq T\}$\ we denote the
natural filtration generated by $\{B_s\}_{0\leq s\leq T}$\ and
augmented by all P-null sets, i.e.,
$${\mathcal{F}}_s=\sigma\{B_r, r\leq s\}\vee \mathcal{N},\ \  s\in [0, T], $$
where $\cal{N}$ is the set of all P-null subsets, and $T > 0$\ a
fixed real time horizon. For any $n\geq 1,$\ $|z|$ denotes the Euclidean norm of $z\in
    {\mathbb{R}}^{n}$. We introduce the following two spaces of processes:
${\cal{S}}^2(0, T; {\mathbb{R}})$\ is the collection of $(\psi_t)_{0\leq t\leq
T}$\ which is a real-valued adapted c\`{a}dl\`{a}g process such that $E[\sup\limits_{0\leq t\leq T}| \psi_{t} |^2]< +\infty$; and ${\cal{H}}^{2}(0,T;{\mathbb{R}}^{n})$\ is the collection of $(\psi_t)_{0\leq t\leq T}$\ which is an ${\mathbb{R}}^{n}$-valued progressively measurable process such that $\parallel\psi\parallel^2_{2}=E[\int^T_0| \psi_t| ^2dt]<+\infty.$

Let $\{A_t, t\geq 0\}$ be a continuous  increasing
${\mathbb{F}}$-progressively measurable scalar process, satisfying
${A}_0=0$ and $E[e^{\mu{A}_T}]<\infty$ for all $\mu>0$. We are given
a final condition $ \xi \in L^{2}(\Omega, {\cal{F}}_{T}, P)$\ such
that $E(e^{\mu A_T}|\xi|^2)<\infty$ for all $\mu>0$, and two
random fields $f: \Omega\times[0, T]\times {\mathbb{R}}\times
{\mathbb{R}}^d\rightarrow {\mathbb{R}}$ and $g: \Omega\times[0,
T]\times {\mathbb{R}}\rightarrow {\mathbb{R}}$ satisfying,
$$\begin{array}{lll}
&{\rm{(i)}}\ \mbox{The processes}\  f(\cdot, y, z)\ \mbox{and}\ g(\cdot, y)\ \mbox{are}\ {\mathbb{F}}\mbox{-progressively measurable and}\\
&\ \ \ \ \ \ \ E[\int_0^Te^{\mu A_t}|f(t,0,0)|^2\,dt]+E[\int_0^Te^{\mu A_t}|g(t,0)|^2dA_t]<\infty,\ \mbox{for all}\ \mu>0;\\
&{\rm{(ii)}}\ \mbox{There is a constant}\ C \ \mbox{such that, for all}\ (t,y,z)\in[0, T]\times{\mathbb{R}}\times {\mathbb{R}}^d,\\
&\ \ \ \ \ \ \ |f(t,y,z)-f(t,y',z')|\leq C(|y-y'|+|z-z'|);\\
&{\rm{(iii)}}\ \mbox{There is a constant}\ C \ \mbox{such that, for all}\ (t,y)\in[0, T]\times{\mathbb{R}},\\
&\ \ \ \ \ \ \ |g(t,y)-g(t,y')|\leq C|y-y'|. \hfill {\rm (H2.1)}
 \end{array}
$$
A solution to the following GBSDE \bel{2.2}Y_t = \xi +
\int_t^Tf(s,Y_s,Z_s)\, ds+\int_t^Tg(s, Y_s)\, dA_s - \int^T_tZ_s\,
dB_s,\quad 0\leq t\leq T,\ee is a pair of ${\mathbb{F}}$-progressively measurable processes $(Y_t, Z_t)_{0\leq t\leq
T}$\ taking values in ${\mathbb{R}}\times {\mathbb{R}}^d$\ which satisfies equation (\ref{2.2}) and
 \be E[\sup\limits_{0\leq t\leq T}|Y_t|^2]+E[\int_0^T|Z_t|^2dt]< \infty,  \ \ 0\leq t\leq T.\ee

From Theorem 1.6 and Proposition 1.1 of~\cite{PZ}, we have the following two lemmas.

 \bl Let (H2.1) be satisfied. Then GBSDE~(\ref{2.2}) has a unique solution $(Y,Z).$\el

\bl Under the assumption (H2.1), we have for any $\mu>0$ \be
 \begin{array}{lll}
 & E[\sup\limits_{0\leq t\leq T}e^{\mu{A}_t}|Y_t|^2+
 \int_0^Te^{\mu{A}_t}|Y_t|^2\, d{A}_t+\int_0^Te^{\mu{A}_t}|Z_t|^2\, dt]\\
                             \leq & C E[e^{\mu{A}_T}|\xi|^2+\int_0^Te^{\mu{A}_t}|f(t,0,0)|^2\, dt
                             +\int_0^Te^{\mu{A}_t}|g(t, 0)|^2\, dA_t]\\
\end{array}\ee for a positive constant $C$, which depends on the Lipschitz
constant of $f$ and $g$, $\mu$, and $T$.\el

Let two sets of data $(\xi, f, g, A)$ and $(\xi', f', g', A')$ satisfy assumption (H2.1). Let $(Y, Z)$ is a solution to GBSDE~(\ref{2.2}) for
data $(\xi, f, g, A)$ and  $(Y', Z')$ for data $(\xi', f', g', A')$. We define
$$(\bar{Y},\ \bar{Z},\ \bar{\xi},\ \bar{f},\ \bar{g},\ \bar{A})=(Y-Y', Z-Z', \xi-\xi', f-f', g-g', A-A').$$ The following two lemmas are borrowed from Proposition 1.2 and Theorem 1.4 of Pardoux and
Zhang~\cite{PZ}, respectively.

\bl For any $\mu>0$, there exists a constant C such that \be
 \begin{array}{lll}
 & E[\sup\limits_{0\leq t\leq T}e^{\mu k_t}|\overline{Y}_t|^2+\int_0^Te^{\mu k_t}|\bar{Z}_t|^2dt]\\
                             \leq &C E[e^{\mu k_T}|\bar{\xi}|^2+\int_0^Te^{\mu
                             k_t}|\bar{f}(t,Y_t,Z_t)|^2dt+\int_0^Te^{\mu k_t}|\bar{g}(t, Y_t)|^2dA'_t
                             +\int_0^Te^{\mu k_t}|g(t, Y_t)|^2d||\bar{A}||_t],\\
   \end{array}
\ee where $k_t:=||\bar{A}||_t+A'_t$, and $||\bar{A}||_t$ is the total
variation of the process $\bar{A}$ on the interval $[0, t]$. \el

For the particular case $A\equiv A',$ we have

\bl(Comparison Theorem) Assume that $\xi\leq\xi',\ f(t,y,z)\leq f'(t,y,z)$, and $g(t,y)\leq g'(t,y),$\ for all $(y, z)\in {\mathbb{R}}\times
{\mathbb{R}}^d,\ dP\times dt, \ \mbox{a.s.}$\ Then $Y_t\leq Y'_t,\
0\leq t\leq T,$\ a.s. \\
Moreover, if $Y_0=Y'_0$, then $Y_t=Y'_t,\ 0 \leq t\leq T$, a.s. In particular, if in addition either $P(\xi<\xi')>0$\ or $f(t,y,z)<f'(t,y,z)$
for any $(y,z)\in {\mathbb{R}}\times {\mathbb{R}}^d$ holds on a set of positive $dt\times dP$\ measure, or $g(t,y)<g'(t,y)$ for any $y\in
{\mathbb{R}}$ holds on a set of positive $dA_t\times dP$\ measure, then $Y_0<Y'_0.$ \el

\section{Formulation of the problem and related DPP}\label{sec3}

\ \ \ \ \ For an admissible control $u(\cdot)\in {\mathcal{U}}$, the
corresponding state process  starting  from $\zeta \in L^2 (\Omega
,{\mathcal{F}}_t, P;\bar{D})$ at the initial time $t$, is governed
by the following reflected SDE:
 \be\label{3.1}
\left\{
  \begin{array}{rcl}
  X_s^{t,\zeta ;u}&=& \zeta+\int_t^sb(r,X_r^{t,\zeta ;u},u_r)dr+\int_t^s\sigma(r,X_r^{t,\zeta
  ;u},u_r)dB_r\\
  && +\int_t^s\nabla\phi(X_r^{t,\zeta ;u})dK_r^{t,\zeta ;u},\ s\in[t,T],\\
  K_s^{t,\zeta ;u}&=& \int_t^sI_{\{X_r^{t,\zeta ;u}\in\partial D\}}dK_r^{t,\zeta ;u},\ K^{t,\zeta ;u}\ \mbox{is
  increasing}.
  \end{array}
  \right.
  \ee
Here, we have made the following assumption on the drift $b:[0,T]\times {\mathbb{R}}^d\times U \rightarrow {\mathbb{R}}^d$
  and the diffusion $\sigma: [0,T]\times {\mathbb{R}}^d\times U \rightarrow {\mathbb{R}}^{d\times d}$:
$$\begin{array}{lll}
&{\rm{(i)}}\ \mbox{For every fixed}\ x\in {\mathbb{R}}^n,\ u\in U,\ b(., x, u)\ \mbox{and}\ \sigma(., x, u)\ \mbox{are continuous in}\ t;\\
&{\rm{(ii)}}\ \mbox{There exists a}\ C>0\ \mbox{such that, for all}\  t\in [0,T],\ x,\ x'\in {\mathbb{R}}^n,\ u \in U, \\
&\ \ \ \ \ |b(t,x,u)-b(t,x',u)|+ |\sigma(t,x,u)-\sigma(t,x',u)|\leq C|x-x'|;\\
&{\rm{(iii)}}\ \mbox{There is some}\ C>0\  \mbox{such that, for all}\  t \in [0, T],\ u\in U\ \mbox{and}\ x\in {\mathbb{R}}^n,\\
&\ \ \ \ \ |b(t,x,u)| +|\sigma (t,x,u)| \leq C(1+|x|).\ \ \ \ \ \ \ \ \ \hfill {\rm (H3.1)}
\end{array}
$$
Therefore, in view of Proposition 5.1 in the Appendix, SDE~(\ref{3.1}) has a unique  strong solution $(X^{t,\zeta ;u},\ K^{t,\zeta ;u})$\ for
any $u(\cdot)\in {\mathcal{U}}$. Moreover, for any $t \in [0,T]$, $u(\cdot)\in {\mathcal{U}},$ and $ \zeta,\ \zeta'\in L^2(\Omega
,{\mathcal{F}}_t,P;\bar{D}),$\ we have
 \be\label{3.3}
\begin{array}{lll}
&E[\sup \limits_{s\in [t,T]}|X^{t,\zeta; u}_s
-X^{t,\zeta';u}_s|^4|{{\mathcal{F}}_t}]\leq C|\zeta -\zeta'|^4, \\
&E[\sup \limits_{s\in [t,T]} |X^{t,\zeta;u}_s|^4|{{\mathcal{F}}_t}]\leq C(1+|\zeta|^4).
\end{array}
\ee Here, the constant $C$ depends only on the Lipschitz
 and the linear growth constants of $b$\ and $\sigma$
with respect to $x$.

We assume that three functions $\Phi: {\mathbb{R}}^d \rightarrow {\mathbb{R}},\ f:[0,T]\times
{\mathbb{R}}^d \times {\mathbb{R}} \times {\mathbb{R}}^d \times U
\rightarrow {\mathbb{R}},\ g:[0,T]\times {\mathbb{R}}^d \times
{\mathbb{R}} \rightarrow {\mathbb{R}}$\ satisfy the following conditions:
$$
\begin{array}{lll}
&{{\rm (i)}}\ \mbox{For every fixed}\ (x, y, z, u)\in {\mathbb{R}}^d \times {\mathbb{R}} \times {\mathbb{R}}^d \times U,\ f(., x, y, z,u)\
\mbox{ is continuous in}\ t ;\ g(\cdot)\in \\
&\ \ \ \ C^{1,2,2}([0,T]\times {\mathbb{R}}^d\times {\mathbb{R}});\ \mbox{and there exists a constant}\  C>0 \ \mbox{such that, for all}\ t\in [0,T],\ x,\\
&\ \ \ \ x'\in {\mathbb{R}}^d,\ y,\ y'\in {\mathbb{R}},\ \ z,\ z'\in {\mathbb{R}}^d,\ u \in U,\\
&\ \  |f(t,x,y,z,u)-f(t,x',y',z',u)|+|g(t,x,y)-g(t,x',y')| \leq C(|x-x'|+|y-y'| +|z-z'|);\\
&{{\rm (ii)}}\ \mbox{There is a constant}\ C>0\ \mbox{such that, for all}\  x,\ x'\in {\mathbb{R}}^d,\\
&\ \ \ \  \ \ |\Phi (x) -\Phi (x')|\leq C|x-x'|;\\
&{{\rm (iii)}}\ \mbox{There exists some}\ C>0\ \mbox{such that, for all}\ 0\leq t \leq T,\  u\in U\ \mbox{and}\  x\in {\mathbb{R}}^n,\\
&\ \ \ \  \ \ |f(t,x,0,0,u)|\leq C(1+|x|).\hfill {\rm (H3.2)}\\
 \end{array}
$$
Then, obviously, $g$\ and $\Phi$ also have the global linear growth condition in $x$: There exists some $C>0$\ such that, for all $0
\leq t \leq T$, and $ x\in {\mathbb{R}}^n $, $|g(t,x,0)|+|\Phi (x)| \leq C(1+|x|).$

For any $u(\cdot) \in {\mathcal{U}}$, and $\zeta \in L^2 (\Omega,{\mathcal{F}}_t,P; \bar{D})$,  the mappings $\xi:= \Phi(X^{t,\zeta; u}_T)$,\ $\widetilde g(s,y):= g(s,X^{t,\zeta; u}_s,y)$\ and $\widetilde f(s,y,z):= f(s,X^{t,\zeta; u}_s,y,z,u_s)$\ satisfy the conditions (H2.1) on the interval $[t, T]$. Therefore, there is a unique solution to the following GBSDE:
      \be\label{3.5}
   \left \{\begin{array}{rcl}
   -dY^{t,\zeta; u}_s & = & \displaystyle f(s,X^{t,\zeta; u}_s, Y^{t,\zeta; u}_s, Z^{t,\zeta; u}_s,
                              u_s)\, ds\\
   && \displaystyle + g(s,X^{t,\zeta; u}_s, Y^{t,\zeta; u}_s)\, dK_s^{t,\zeta; u}-Z^{t,\zeta; u}_s\, dB_s,\\
      Y^{t,\zeta; u}_T  & = &\displaystyle  \Phi (X^{t,\zeta; u}_T),
   \end{array}\right.
   \ee
where $(X^{t,\zeta; u},\ K^{t,\zeta; u})$ solves the reflected SDE~(\ref{3.1}).

 Moreover, similar to Proposition 5.2, there exists some constant $C>0$\ such that, for all $0 \leq t \leq
T,\ \zeta,\ \zeta' \in L^2(\Omega , {\mathcal{F}}_t,P;\bar{D}),\ u(\cdot) \in {\mathcal{U}}, $\ $P$-a.s.,
 \be\label{3.6}
\begin{array}{ll}
 {\rm(i)} & |Y^{t,\zeta; u}_t -Y^{t,\zeta'; u}_t| \leq C(|\zeta -\zeta'|+|\zeta -\zeta'|^{\frac{1}{2}}); \\
 {\rm(ii)} & |Y^{t,\zeta; u}_t| \leq C (1+|\zeta|). \\
\end{array}
\ee

 We now introduce the following definitions about admissible controls.

\noindent\begin{definition}\label{d3.1} An admissible control process $u=\{u_r, r\in [t, s]\}$ on $[t, s]$ (with $s\in (t,T]$) is an
${\mathcal{F}}_r$-progressively measurable process taking values in U. The set of all admissible controls  on $[t, s]$ is denoted by\
${\mathcal{U}}_{t, s}.$\  We identify two processes $u$\ and $\bar{u}$ in ${\mathcal{U}}_{t, s}$ and write $u\equiv \bar{u}\ \mbox{on}\ [t,
s],$\ if $P\{u=\bar{u}\ \mbox{a.e. in}\ [t, s]\}=1.$\
\end{definition}

 For any $u(\cdot)\in {\mathcal{U}}_{t,T}$,   the value of
 the associated cost functional is given by
 \be\label{3.7}
J(t, x; u):= Y^{t, x; u}_t,\ (t, x)\in [0, T]\times \bar{D},\ee
where the process $Y^{t, x; u}$ is defined by GBSDE~(\ref{3.5}).

From Theorem \ref{th6.1}, we have
 \be\label{3.8} J(t, \zeta; u) = Y^{t,\zeta; u}_t,\quad (t,\zeta)\in[0,
T]\times L^2 (\Omega ,{\mathcal{F}}_t ,P; \bar{D}).\ee We define the
value function of our stochastic control problem as follows:
\be\label{3.9} W(t,x):= \esssup_{u \in {\mathcal{U}}_{t,T}}J(t,x;
u), \quad (t,x)\in [0,T]\times \bar{D}. \ee

Under assumptions (H3.1) and (H3.2), the value function $W$ is well-defined on $[0,T]\times D$, and its values at time $t$ are bounded and
${\mathcal{F}}_{t}$-measurable random variables. In fact, they are all deterministic. We have

 \begin{proposition}\label{p3.1} For any $(t, x)\in [0, T]\times \bar{D}$, we have $W(t,x)=E[W(t,x)]$, P-a.s. Let $W(t,x)$\ equal to its deterministic version $E[W(t,x)]$.  Then $W:[0, T]\times \bar{D}\to \mathbb{R}$\ is a deterministic function.\end{proposition}

The proof is an adaptation of relevant arguments of Buckdahn and Li~\cite{BL}. For the readers' convenience we give it in the Section 5.2 of Appendix.

 As an immediate result of (\ref{3.6}) and (\ref{3.9}), the value function $W$ has the following property .

\begin{lemma}\label{l3.2}There exists a constant $C>0$\ such that, for all $ (t, x, x')\in [0,T]\times \bar{D}\times \bar{D}$, \be\label{3.12}
\begin{array}{llll}
&{\rm(i)} & |W(t,x)-W(t,x')| \leq C[|x-x'|+|x-x'|^{\frac{1}{2}}];  \\
&{\rm(ii)} & |W(t,x)| \leq C(1+|x|).
\end{array}
\ee \end{lemma}

We now study the (generalized) DPP for our stochastic control
problem~(\ref{3.1}), (\ref{3.5}), and (\ref{3.9}). For this we
have to define the family of (backward) semigroups related with
GBSDE~(\ref{3.5}). Peng~\cite{Pe1} first introduced the notion of backward stochastic semigroups to study the DPP for the optimal
stochastic control of SDEs. In what follows, it is adapted to the optimal control problem of stochastic differential systems  {\it reflected}  in a domain.

 Given the initial data $(t,x)$, a positive number $\delta\leq T-t$, an admissible control $u(\cdot) \in {\mathcal{U}}_{t, t+\delta}$, and a random variable $\eta \in L^2 (\Omega, {\mathcal{F}}_{t+\delta},P; {\mathbb{R}})$, we define \be\label{3.13}
G^{t, x; u}_{s,t+\delta} [\eta]:= \tilde{Y}_s^{t,x; u},\ \hskip0.5cm
s\in[t, t+\delta], \ee where $(\tilde{Y}_s^{t,x;u},
\tilde{Z}_s^{t,x;u})_{t\leq s \leq t+\delta}$ is the solution of the
following GBSDE on the time interval $[t, t+\delta]$:
$$
\left \{\begin{array}{rcl}
 -d\tilde{Y}_s^{t,x;u} & = & f(s,X^{t,x;u}_s ,\tilde{Y}_s^{t,x;u}, \tilde{Z}_s^{t,x;u},
                           u_{s})\, ds+ g(s,X^{t,x;u}_s ,\tilde{Y}_s^{t,x;u})\, dK_s^{t,x;u}\\
  && -\tilde{Z}_s^{t,x; u} dB_s , \quad s\in [t,t+\delta];\\
 \tilde{Y}_{t+\delta}^{t,x; u} & =& \eta ,
\end{array}\right.
$$
and $(X^{t,x;u}, K^{t,x;u})$\ is the solution of reflected SDE~(\ref{3.1}). Then, obviously, for the solution $(Y^{t,x;u}, Z^{t,x;u})$\ of
GBSDE~(\ref{3.5}), we have \be\label{3.14} G^{t,x;u}_{t,T} [\Phi (X^{t,x; u}_T)] =G^{t,x;u}_{t,t+\delta} [Y^{t,x;u}_{t+\delta}]. \ee
Furthermore,
$$
\begin{array}{rcl}
 J(t,x;u)=Y_t^{t,x;u}=G^{t,x;u}_{t,T} [\Phi (X^{t,x; u}_T)]
  =G^{t,x;u}_{t,t+\delta} [Y^{t,x;u}_{t+\delta}]=G^{t,x;u}_{t,t+\delta} [J(t+\delta,X^{t,x;u}_{t+\delta};u)].
\end{array}
$$

\begin{remark}\label{r3.3} If both $f$ and $g$ do not depend on  $(y, z)$, we have
$$G^{t,x;u}_{s,t+\delta}[\eta]=E[\eta + \int_s^{t+\delta}f(r,X^{t,x;u}_r,u_{r})dr + \int_s^{t+\delta}g(r,X^{t,x;u}_r)dK^{t,x;u}_r|{\cal{F}}_s],\ \ s\in [t,
t+\delta].$$
\end{remark}
\begin{theorem}\label{th3.1} Under
assumptions {\rm (H3.1)} and {\rm (H3.2)}, the value function $W$
satisfies the following DPP: For any $0\leq t<t+\delta \leq T,\ x\in \bar{D},$
 \be\label{3.15}
W(t,x) =\esssup_{u \in {\mathcal{U}}_{t,
t+\delta}}G^{t,x;u}_{t,t+\delta} [W(t+\delta,
X^{t,x;u}_{t+\delta})].
 \ee
  \end{theorem}

The proof is similar to [4]. For the readers' convenience we give it in Section 5.2.

Lemma \ref{l3.2} shows that the value function $W(t,x)$ is continuous in $x$, uniformly in $t$. From Theorem \ref{th3.1} we can
get the continuity of $W(t,x)$\ in $t$.

 \begin{theorem}\label{th3.2} Let assumptions {\rm (H3.1)} and {\rm (H3.2)} be satisfied.
 Then the value function $W(t,x)$ is continuous in $t$.
  \end{theorem}

\noindent{\it Proof}. Let $(t, x)\in [0,T]\times \bar{D}$\ and $\delta \in
(0, T-t]$. We want to prove that $W$ is continuous in $t$.
For this we notice that from (\ref{3.28}), for an arbitrarily small $\varepsilon>0,$
\be\label{3.31} I^1_\delta +I^2_\delta \leq W(t,x)-W(t+\delta ,x)
\leq I^1_\delta +I^2_\delta +C\varepsilon, \ee
 where
$$
\begin{array}{lll}
I^1_\delta & := & G^{t,x; u^{\varepsilon}}_{t,t+\delta}[W(t+\delta,
X^{t,x; u^{\varepsilon}}_{t+\delta})]
                   -G^{t,x;u^{\varepsilon}}_{t,t+\delta} [W(t+\delta,x)], \\
I^2_\delta & := & G^{t,x; u^{\varepsilon}}_{t,t+\delta}[W(t+\delta,x)] -W(t+\delta ,x),
\end{array}
$$
for $u^{\varepsilon} \in {\cal{U}}_{t, t+\delta}$\ such that (\ref{3.28}) holds. From Lemma 2.3 and the estimate (\ref{3.12}) we get that, for
some constant $C$ which does not depend on the controls $u^{\varepsilon}$,
$$
\begin{array}{rcl}
|I^1_\delta | &\leq& (CE[|W(t+\delta ,X^{t,x;
u^{\varepsilon}}_{t+\delta})
                  -W(t+\delta ,x)|^2|{{\mathcal{F}}_t}])^{\frac{1}{2}}\\
              & \leq& (CE[|X^{t,x;u^{\varepsilon}}_{t+\delta} -x|^2+|X^{t,x;u^{\varepsilon}}_{t+\delta} -x||{{\mathcal{F}}_t}])^{\frac{1}{2}},
\end{array}
$$
and since $E[|X^{t,x;u^{\varepsilon}}_{t+\delta}
-x|^2|{{\mathcal{F}}_t}] \leq C\delta $ (refer to (\ref{5.15}) in Appendix) we get that $|I^1_\delta| \leq C\delta^{\frac{1}{4}}$. From the definition of $G^{t,x; u^{\varepsilon}}_{t,t+\delta}[\cdot]$\ (see (\ref{3.13})),
$$\begin{array}{rcl}
I^2_\delta  & = & E[W(t+\delta ,x) +\int^{t+\delta}_t
f(s,X^{t,x;u^{\varepsilon}}_s,\tilde{Y}^{t,x;u^{\varepsilon}}_s,
\tilde{Z}^{t,x;u^{\varepsilon}}_s, u^{\varepsilon}_s) ds   \\
 &         & +\int^{t+\delta}_t g(s,X^{t,x;u^{\varepsilon}}_s,\tilde{Y}^{t,x;u^{\varepsilon}}_s)dK^{t,x;u^{\varepsilon}}_s
-\int^{t+\delta}_t \tilde{Z}^{t,x; u^{\varepsilon}}_s dB_s|{{\mathcal{F}}_t}] -W(t+\delta ,x)  \\
 &    =  &  E[\int^{t+\delta}_t f(s,X^{t,x; u^{\varepsilon}}_s,\tilde{Y}^{t,x; u^{\varepsilon}}_s,
 \tilde{Z}^{t,x; u^{\varepsilon}}_s, u^{\varepsilon}_s)ds +\int^{t+\delta}_t g(s,X^{t,x;u^{\varepsilon}}_s,\tilde{Y}^{t,x;u^{\varepsilon}}_s)
dK^{t,x;u^{\varepsilon}}_s|{{\mathcal{F}}_t}].
\end{array}$$
From the Schwartz inequality, Propositions 5.2 and 5.3 in Appendix and (\ref{3.3}), we then get
$$\begin{array}{rcl}
|I^2_\delta | &\leq& \delta^{\frac{1}{2}} E[\int^{t+\delta}_t |f(s,X^{t,x;u^{\varepsilon}}_s,
     \tilde{Y}^{t,x;u^{\varepsilon}}_s,\tilde{Z}^{t,x;u^{\varepsilon}}_s,
u^{\varepsilon}_s)|^2ds|{{\mathcal{F}}_t}]^{\frac{1}{2}}  \\
& &+E[K_{t+\delta}^{t,x;u^{\varepsilon}}|{{\mathcal{F}}_t}]^{\frac{1}{2}}
     E[\int^{t+\delta}_t |g(s,X^{t,x;u^{\varepsilon}}_s,
     \tilde{Y}^{t,x;u^{\varepsilon}}_s)|^2dK_s^{t,x;u^{\varepsilon}}|{{\mathcal{F}}_t}]^{\frac{1}{2}}  \\
 &\leq& \delta^{\frac{1}{2}}E[\int^{t+\delta}_t
(|f(s,X^{t,x;u^{\varepsilon}}_s,0,0,u^{\varepsilon}_s)|+C|
\tilde{Y}^{t,x;u^{\varepsilon}}_s|+ C|\tilde{Z}^{t,x;u^{\varepsilon}}_s|)^2ds|{{\mathcal{F}}_t}]^{\frac{1}{2}}\\
& &+E[K_{t+\delta}^{t,x;u^{\varepsilon}}|{{\mathcal{F}}_t}]^{\frac{1}{2}}E[\int^{t+\delta}_t
(|g(s,X^{t,x;u^{\varepsilon}}_s,0)|+C|\tilde{Y}^{t,x;u^{\varepsilon}}_s|)^2dK_s^{t,x;u^{\varepsilon}}|{{\mathcal{F}}_t}]^{\frac{1}{2}}\\
 &\leq& C\delta^{\frac{1}{2}}E[\int^{t+\delta}_t(|1+|X^{t,x;u^{\varepsilon}}_s|+| \tilde{Y}^{t,x;u^{\varepsilon}}_s|
+|\tilde{Z}^{t,x;u^{\varepsilon}}_s|)^2ds|{{\mathcal{F}}_t}]^{\frac{1}{2}}\\
& &+CE[K_{t+\delta}^{t,x;u^{\varepsilon}}|{{\mathcal{F}}_t}]^{\frac{1}{2}}E[\int^{t+\delta}_t(|1+|X^{t,x;u^{\varepsilon}}_s|+| \tilde{Y}^{t,x;u^{\varepsilon}}_s|)^2dK_{s}^{t,x;u^{\varepsilon}}|{{\mathcal{F}}_t}]^{\frac{1}{2}}\\
 &\leq& C\delta^{\frac{1}{2}}+C(E[|K_{t+\delta}^{t,x;u^{\varepsilon}}|^2|{{\mathcal{F}}_t}])^2
 \leq C\delta^{\frac{1}{2}}.
\end{array}$$
Then, from (\ref{3.31}), $|W(t,x)-W(t+\delta ,x)| \leq C \delta^{\frac{1}{4}}+C\delta^{\frac{1}{2}}+C\varepsilon,$\ and letting
$\varepsilon \downarrow 0$\ we get $W(t,x)$\ is continuous in $t$.
The proof is complete. \endproof

\section{Viscosity solutions of related HJB equations}\label{sec4}

\ \ \ \ \ We consider the following PDE: \be\label{4.1}
 \left \{\begin{array}{ll}
 &\!\!\!\!\! \frac{\partial }{\partial t} W(t,x) +  H(t, x, W, DW,
 D^2W)=0, \quad\quad \hfill (t,x)\in [0,T)\times {D} ,  \\
 &\!\!\!\!\! \frac{\partial }{\partial n} W(t,x)+g(t,x,W(t,x))=0, \hfill 0\leq t< T,\ x\in \partial{D};\\
 &\!\!\!\!\!   W(T,x) =\Phi (x), \hfill   x \in \bar{D},
 \end{array}\right.
\ee  where at a point $x\in \partial D$, $\frac{\partial }{\partial
n}=\sum_{i=1}^{d}\frac{\partial }{\partial x_i}\phi(x)\frac{\partial
}{\partial x_i}$,\ and the Hamiltonian $H$ is defined by
$$ H(t, x, y, p, A):= \sup_{u
\in U}\{\frac{1}{2}{\text{tr}}(\sigma\sigma^{T}(t, x,
 u)A)+ \langle p, b(t, x, u)\rangle+ f(t, x, y, p\sigma,
u)\},$$  where $ (t,x,y,p,A)\in [0, T]\times {\mathbb{R}}^n\times {\mathbb{R}}\times {\mathbb{R}}^d\times {\mathbf{S}}^d$
with ${\mathbf{S}}^d$ being the set of all $d\times d$ symmetric matrices.

 In this section we shall prove that the value function $W$ defined by (\ref{3.9}) is the unique
 viscosity solution of (\ref{4.1}). The interested reader is referred to Crandall, Ishii, and
Lions~\cite{CIL} for a detailed introduction to viscosity solutions. Let $C^3_{l, b}([0,T]\times \bar{D})$ be the set of the
real-valued functions that are continuously differentiable up to the
third order and whose derivatives of order from 1 to 3 are bounded.

\begin{definition}\label{d4.1} A real-valued
continuous function $W\in C([0,T]\times \bar{D} )$ is called \\
  {\rm(i)} a viscosity subsolution of {\rm (\ref{4.1})} if $W(T,x) \leq \Phi (x),\ \mbox{for all}\ x \in
  \bar{D}$, and if for all functions $\varphi \in C^3_{l, b}([0,T]\times
  \bar{D})$ and $(t,x) \in [0,T) \times \bar{D}$ such that $W-\varphi $\ attains its
  local maximum at $(t, x)$:
     $$\begin{array}{lll}
    & \frac{\partial \varphi}{\partial t} (t,x)+  H(t, x, W, D\varphi, D^2\varphi) \geq 0, \  \mbox{ \rm if } x\in
    D;\\
& \max \{\frac{\partial \varphi}{\partial t} (t,x)
+  H(t, x, W, D\varphi, D^2\varphi), \ \ \frac{\partial
\varphi}{\partial n} (t,x) + g(t, x, W)\} \geq
0, \mbox{ \rm if } x\in\partial{D};\\
     \end{array}$$
{\rm(ii)} a viscosity supersolution of {\rm (\ref{4.1})} if $W(T,x)
\geq \Phi (x),~\mbox{for all}\ x \in
  \bar{D}$, and if for all functions $\varphi \in C^3_{l, b}([0,T]\times
 \bar{D})$ and $(t,x) \in [0,T) \times \bar{D}$ such that $W-\varphi $\ attains its
  local minimum at $(t, x)$:
$$\begin{array}{lll}
    & \frac{\partial \varphi}{\partial t} (t,x) +  H(t, x, W, D\varphi, D^2\varphi) \leq 0,\ \  \mbox{if}\ \ x\in
    D;\\
& \min \{\frac{\partial \varphi}{\partial t} (t,x) +  H(t, x, W,
D\varphi, D^2\varphi), \ \ \frac{\partial \varphi}{\partial n} (t,x)
+ g(t, x, W) \} \leq
0,\ \mbox{if}\ x\in\partial{D};\\
     \end{array}$$
{\rm(iii)} a viscosity solution of {\rm (\ref{4.1})} if it is both a
viscosity sub- and a supersolution of
     {\rm (\ref{4.1})}.\end{definition}

For simplicity of notations, we define for $\varphi \in C^3_{l, b}([0,T] \times \bar{D})$,
\begin{equation}\begin{array}{rcl}\label{4.3}
     F(s,x,y,z,u)&£º=& \frac{\partial }{\partial s}\varphi (s,x) + \frac{1}{2}{\rm tr}(\sigma\sigma^{T}(s,x, u)D^2\varphi)+ D\varphi.b(s, x, u) \\
       & & + f(s, x, y+\varphi (s,x), z+ D\varphi (s,x).\sigma(s,x,u),u), \\
     G(s,x,y)&£º=& \frac{\partial }{\partial n}\varphi (s,x) +
      g(s, x, y+\varphi (s,x)), \end{array}
     \end{equation}
for $(s, x, y, z, u)\in [0,T] \times \bar{D} \times {\mathbb{R}} \times {\mathbb{R}}^d \times U$.

\begin{proposition}\label{pr4.1}Under the assumptions {\rm (H3.1)} and {\rm (H3.2)} the value function $W$
is a viscosity subsolution to {\rm (\ref{4.1})}. \end{proposition}

\noindent{\it Proof}. Obviously, $W(T,x)=\Phi (x),\ x\in \bar{D}$. Suppose
that $\varphi \in C^3_{l,b} ([0,T] \times \bar{D })$ and $(t,x)\in
[0, T)\times \bar{D }$ is such that $W-\varphi$\ attains its maximum
at $(t,x)$. Without loss of generality, we assume that $\varphi
(t,x)=W(t,x)$.

We first consider the case $x\in D$. We shall prove that $$\sup_{u\in U}F(t,x,0,0,u)\geq 0. $$ If this is not true, then there exists some $\theta>0$ such that
\be\label{4.16} F_0(t,x):=\sup_{u\in U}
F(t,x,0,0,u)\leq-\theta<0.\ee Therefore, $ F(t,x,0,0,u)\leq -\theta,\ \mbox{for all}\
u\in U.$

Since $F_0$ is continuous at $(t,x)$, we can choose $\bar{\alpha}\in
(0, T-t]$ such that
\begin{eqnarray}
&O_{\bar{\alpha}}(x):=\{y:|y-x|\leq \bar{\alpha}\}\subset D, \\
&F(s,y,0,0,u)\leq-\frac{1}{2}\theta,\ \mbox{ \rm for all }(s,y,u)\in
[t, t+\bar{\alpha}]\times O_{\bar{\alpha}}(x) \times U.\label{F1}
\end{eqnarray}
For any $\alpha\in (0, \bar{\alpha}]$, we consider the following BSDE:
    \be\label{4.4} \left \{\begin{array}{rcl}
      -dY^{1,u}_s &=&\displaystyle F(s,X^{t,x;u}_s,
      Y^{1,u}_s,Z^{1,u}_s,u_s)\, ds+G(s,X^{t,x;u}_s,
      Y^{1,u}_s)\, dK_s^{t,x;u}\\
      &&\displaystyle  -Z^{1,u}_s dB_s, \quad s\in [t,t+\alpha]; \\
     Y^{1,u}_{t+\alpha}&=&0,
     \end{array}\right.\ee
     where the pair of processes $(X^{t,x,u},K^{t,x,u})$ are given by
     $(\ref{3.1})$\ and $u(\cdot) \in {\mathcal{U}}_{t, t+\alpha}$.
It is not hard to check that $F(s,X^{t,x;u}_s, y,z,u_s)$
and $ G(s,X^{t,x;u}_s, y)$ satisfy (H2.1). Thus, due to
Lemma 2.1, GBSDE~(\ref{4.4}) has a unique solution. We have the
following observation.

\begin{lemma}\label{l4.1} For every $s\in [t,t+\alpha]$, we have the following
relationship:
    \be\label{4.5}
     Y^{1,u}_s = G^{t,x;u}_{s,t+\alpha} [\varphi (t+\alpha ,X^{t,x;u}_{t+\alpha})]
                -\varphi (s,X^{t,x;u}_s), \hskip 0.5cm
              \mbox{{\it P}-a.s.}\ee
 \end{lemma}

\noindent {\it Proof}. We recall that $G^{t,x;u}_{s,t+\alpha} [\varphi
(t+\alpha, X^{t,x;u}_{t+\alpha})]$ is defined by the
solution of the GBSDE
     $$
     \left \{\begin{array}{rcl}
     -dY^{u}_s &=& f(s,X^{t,x;u}_s, Y^{u}_s,Z^{u}_s,u_s)\, ds+g(s,X^{t,x;u}_s,
     Y^{u}_s)\, dK_s^{t,x;u}\\
     && -Z^{u}_s\, dB_s , \hskip 0.2cm s\in [t,t+\alpha]; \\
     Y^{u}_{t+\alpha}&=& \varphi (t+\alpha
     ,X^{t,x;u}_{t+\alpha}),
     \end{array}\right.
     $$
with the following formula:
     \be\label{4.6}
     G^{t,x;u}_{s,t+\alpha} [\varphi (t+\alpha ,X^{t,x;u}_{t+\alpha})] =Y^{u}_s, \hskip 0.5cm
     s\in [t,t+\alpha],  \ee
(see (\ref{3.13})). Hence, we only need to show that
$Y^{u}_s-\varphi (s,X^{t,x;u}_s)\equiv Y^{1,u}_s$ for $s\in
[t,t+\alpha]$. This can be verified directly by applying It\^{o}'s
formula to $\varphi (s,X^{t,x;u}_s)$. Indeed, the stochastic differentials of $Y^{u}_s -\varphi (s,X^{t,x;u}_s)$ and
$Y^{1,u}_s$\ equal, and with the same terminal condition
$Y^{u}_{t+\alpha} - \varphi (t+\alpha ,X^{t,x;u}_{t+\alpha}) =0 =
Y^{1,u}_{t+\alpha}.$  \endproof
\br For $x\in \partial D$\ Lemma 4.1 still holds.
\er
On the other hand, from the DPP (see Theorem \ref{th3.1}), for
every $\alpha$,
     $$
     \varphi (t,x) =W(t,x) =\esssup_{u \in
{\mathcal{U}}_{t, t+\alpha}}G^{t,x;u}_{t,t+\alpha} [W(t+\alpha,
X^{t,x;u}_{t+\alpha})],
     $$
  and from $W\leq\varphi$\ and the monotonicity property of $G^{t,x;u}_{t,t+\delta}[\cdot]$\ (see Lemma 2.4)\ we get
     $$
     \esssup_{u \in
{\mathcal{U}}_{t, t+\alpha}} \{G^{t,x;u}_{t,t+\delta}
[\varphi(t+\alpha, X^{t,x;u}_{t+\alpha})] -\varphi (t,x)\}\geq 0,\
\mbox{{\it P}-a.s.}
     $$
 Thus, from Lemma \ref{l4.1}, we have $  \esssup_{u \in
{\mathcal{U}}_{t, t+\alpha}} Y^{1,u}_t \geq 0,\quad \mbox{{\it P}-a.s.}$

 Hence, for arbitrary $\varepsilon>0$, similar to that of inequality (\ref{3.28}),  there is $u^\varepsilon\in
{\mathcal{U}}_{t, t+\alpha}$\ such that \be
Y^{1,u^\varepsilon}_t\geq -\varepsilon\alpha, \quad \mbox{{\it P}-a.s.}\label{4.7}\ee

\br Similarly, (\ref{4.7}) is still true for $x\in \partial D$.
\er

For $u^\varepsilon\in {\mathcal{U}}_{t, t+\alpha}$\ we define $\tau=\inf\{s\geq t: |X_s^{t,x;u^\varepsilon}-x|\geq \bar{\alpha}\}\wedge
(t+\alpha).$\ Consequently, on $[t, \tau]$\ the process $(K^{t,x;u})$\ is zero and, hence
$$Y_s^{1;u^\varepsilon}=Y_\tau^{1;u^\varepsilon}+\int_s^\tau F(r,X_r^{t,x;u^\varepsilon},Y_r^{1;u^\varepsilon},
Z_r^{1;u^\varepsilon}, u^\varepsilon_r)dr-\int_s^\tau
Z_r^{1;u^\varepsilon}dB_r. $$ We consider the following two BSDEs:
\be\label{4.10} \left \{\begin{array}{rcl}
      -dY^2_s &=& (C^*(\vert Y^2_s\vert+\vert Z^2_s\vert)-\frac{1}{2}\theta)\, ds-Z^2_s\, dB_s, \\
     Y^2_{t+\alpha}&=&0,
     \end{array}\right.
     \ee
whose unique solution is given by \be\label{4.11}
Y^2_s=-\frac{\theta}{2C^*}(1-e^{C^*(s-(t+\alpha))}),\
Z^2_s=0, \quad s\in[t,t+\alpha],\ee
and \be \label{4.12}\left
\{\begin{array}{rcl}
      -dY^3_s &=& (C^*(\vert Y^3_s\vert+\vert Z^3_s\vert)-\frac{1}{2}\theta)\, ds-Z^3_s\, dB_s, \quad s\in[t,\tau];\\
     Y^3_{\tau}&=& Y_\tau^{1;u^\varepsilon}.
     \end{array}\right.
     \ee
Here, $C^*$ is the Lipschitz constant of $F$\ with respect to $y,\ z$; also the Lipschitz constant of $G$\ with respect to $y$, in order to be different from the constant $C$ which may vary from lines to lines.
We have the following lemma.

\bl We have $ Y_t^{1,u^\varepsilon}
 \leq Y_t^3$ and $|Y_t^2-Y_t^3|\leq C\alpha^{\frac{3}{2}},
 \mbox{P-a.s.}$ Here $C>0$\ is  independent of both the control $u$ and
$\alpha$.\el

\noindent{\it Proof}. (1) We observe from (\ref{F1}) and the definition of $\tau$
that, for all $(s,y,z,u)\in [t,\tau]\times {\mathbb{R}}\times
{\mathbb{R}}^d\times U$,
$$  \begin{array}{llll}
F(s,X_s^{t,x;u^\varepsilon},y,z, u^\varepsilon) &\leq & C^*(\vert y\vert+\vert
z\vert)+F(s,X_s^{t,x;u^\varepsilon},0,0,u^\varepsilon)\\
& \leq & C^*(\vert y\vert+\vert z\vert)-\frac{1}{2}\theta.
       \end{array}$$
Consequently, from Lemma 2.2 in \cite{BL} (the comparison result for BSDEs) we
have that $$Y_s^{1,u^\varepsilon}\leq Y^3_s,\ s\in[t,\tau],\ \mbox{P-a.s.,}$$
 where $Y^3$ is defined by BSDE~(\ref{4.12}).

(2) From the equation (\ref{4.4}), Proposition 5.1 and Proposition 5.2 in the
Appendix, we have
$$|Y_\tau^{1;u^\varepsilon}|\leq C(t+\alpha-\tau)^{\frac{1}{2}}
+C(E[(K_{t+\alpha}^{t,x;u^\varepsilon}-K_{\tau}^{t,x;u^\varepsilon})^2|{\cal{F}}_\tau])^{\frac{1}{2}},$$ where $C$ is independent of controls, and
$K_{t+\alpha}^{t,x;u^\varepsilon}-K_{\tau}^{t,x;u^\varepsilon}=K_{t+\alpha}^{\tau,X_\tau^{t,x;u^\varepsilon};u^\varepsilon}$
by means of the uniqueness of solution of reflected SDE~(\ref{3.1}).
Therefore, we have
$$E[|Y_{\tau}^{1;u^\varepsilon}|^2|{\cal{F}}_t]\leq
CE[(t+\alpha-\tau)|{\cal{F}}_t] +CE[|K_{t+\alpha}^{\tau,X_\tau^{t,x;u^\varepsilon};u^\varepsilon}|^{2}|{\cal{F}}_t].$$ From Proposition 5.3 in Appendix, we have
\be\label{4.13-1}E[|K_{t+\alpha}^{\tau,X_\tau^{t,x;u^\varepsilon};u^\varepsilon}|^{\, 2}|{\cal{F}}_t]\leq
C(E[(t+\alpha-\tau)^2|{\cal{F}}_t])^{\frac{1}{2}}.\ee Therefore, we get
\be\label{4.14-2}E[|Y_{\tau}^{1;u^\varepsilon}|^2|{\cal{F}}_t]\leq
C(E[(t+\alpha-\tau)^2|{\cal{F}}_t])^{\frac{1}{2}}.\ee On the other hand, we
consider the following SDE: \be\label{4.13}\begin{array}{rcl}
&d\bar{X}_s^{t,x;u^\varepsilon}=b(s,
\bar{X}_s^{t,x;u^\varepsilon},u_s^{\varepsilon})\, ds+\sigma(s,
\bar{X}_s^{t,x;u^\varepsilon},u_s^{\varepsilon})\, dB_s,\ s\geq t;
&\bar{X}_t^{t,x;u^\varepsilon}=x.  \end{array}\ee Then we know on
$[t,\tau]$, P-a.s., $ {X}^{t,x;u^\varepsilon}=\bar{X}^{t,x;u^\varepsilon}.$\ For $\bar{X}^{t,x;u^\varepsilon}$ we have the classical estimate
$$E[\sup_{t\leq s\leq t+\alpha}|\bar{X}^{t,x;u^\varepsilon}_s-x|^8|{\cal{F}}_t]\leq C\alpha^4,\ \mbox{P-a.s.}$$
Therefore, we have
 \be\label{4.14}
P\{\tau<t+\alpha|{\cal{F}}_t\}\leq P\{\sup_{s\in [t,
t+\alpha]}|\bar{X}_s^{t,x;u^\varepsilon}-x|\geq\bar{\alpha}|{\cal{F}}_t\}\leq \frac{C}{\bar{\alpha}^8}\alpha^4.\ee
Hence, \be\label{4.15} E[|Y_{\tau}^{1;u^\varepsilon}|^2|{\cal{F}}_t]\leq
C\alpha (P\{\tau<t+\alpha|{\cal{F}}_t\})^{\frac{1}{2}} \leq
\frac{C}{\bar{\alpha}^4}\alpha^3.\ee Furthermore, from Lemma 2.3 in \cite{BL},
\be\label{4.16}\begin{array}{rcl} |Y_t^2-Y_t^3|&\leq&
C(E[|Y_{\tau}^{2}-Y_{\tau}^{3}|^2|{\cal{F}}_t])^{\frac{1}{2}}\leq C(E[|Y_{\tau}^{2}|^2|{\cal{F}}_t])^{\frac{1}{2}}
+C(E[|Y_{\tau}^{3}|^2|{\cal{F}}_t])^{\frac{1}{2}}\\
&\leq& C
\frac{\theta}{2}(1-e^{-C^*\alpha})(P\{\tau<t+\alpha|{\cal{F}}_t\})^{\frac{1}{2}}
+C(E[|Y_{\tau}^{1;u^\varepsilon}|^2|{\cal{F}}_t])^{\frac{1}{2}}\\
&\leq&
C\frac{\theta}{2}(1-e^{-C^*\alpha})\frac{1}{\bar{\alpha}^4}\alpha^2+\frac{C}{\bar{\alpha}^2}\alpha^{\frac{3}{2}}\leq C\alpha^{\frac{3}{2}},\\
\end{array}
\ee for any $\alpha\in (0,\bar{\alpha}].$
\bigskip

\noindent\textbf{Proof of Proposition 4.1 (sequel)}.

By combining (\ref{4.7}) with Lemma 4.2 we then obtain
$$- \varepsilon\alpha\leq Y^{1,u^\varepsilon}_t\leq Y_t^3\leq Y_t^2+|Y_t^2-Y_t^3|\leq Y_t^2+C\alpha^{\frac{3}{2}},\ \mbox{P-a.s.}$$
i.e.,\ \ \ \ $- \varepsilon\alpha\leq Y^{1,u^\varepsilon}_t\leq - \frac{\theta }{2C^*}(1-e^{-
C^*\alpha})+C\alpha^{\frac{3}{2}},\ \mbox{P-a.s.}$\ \ \ Therefore, $$- \varepsilon\leq - \frac{\theta }{2 C^*}\frac{1-e^{-
C^*\alpha}}{\alpha}+C\alpha^{\frac{1}{2}}.$$ Letting $\alpha\to 0+$
and $\varepsilon\to 0+$,  we get $0\leq - \frac{\theta }{2}$, which
contradicts our assumption that $\theta>0$. Therefore, we have\ \ $
     \sup_{u\in U} F(t,x,0,0,u) \geq 0, $\ \  which implies by the definition of $F$  that $$\frac{\partial \varphi}{\partial t} (t,x) +  H(t, x,
W, D\varphi, D^2\varphi) \geq 0,\ \  \mbox{if}\ \ x\in
    D.$$

We now consider the case $x\in \partial D$. We must prove that
$$\max\{\frac{\partial \varphi}{\partial t} (t,x) +  H(t, x,
\varphi, D\varphi, D^2\varphi), \ \ \frac{\partial \varphi}{\partial
n} (t,x) + g(t, x, \varphi)\} \geq 0$$ If this is not true,
then there exists some $\theta>0$ such that \be\label{4.17}
\sup_{u\in
U} F(t,x,0,0,u)\leq-\theta<0,\quad
G(t,x,0)\leq-\theta<0,\ee therefore,\ \ $ F(t,x,0,0,u)\leq -\theta\ \  \mbox{for all}\
u\in U;$\ \ and $ G(t,x,0)\leq -\theta\ \ \mbox{for all}\ \ u\in U.$

Choose $\bar{\alpha}\in (0, T-t]$ such that \begin{eqnarray}
&F(s,y,0,0,u)\leq-\frac{1}{2}\theta,\label{F2} \\
&G(s,y,0)\leq-\frac{1}{2}\theta,\  \mbox{for all}\ u\in U, t\leq
s\leq t+\bar{\alpha},\ |y-x|\leq \bar{\alpha}.\label{G1}
\end{eqnarray}

Now we fix $\bar{\alpha}$, and we consider any $\alpha\in
(0, \bar{\alpha}].$\ Similarly, we consider GBSDE (4.6) with $x\in \partial D$, then we also can get (4.7) and (4.9).
For $u^\varepsilon\in {\mathcal{U}}_{t,
t+\alpha}$\ in (\ref{4.7}) we define
$$\tau=\inf\{s\geq t: |X_s^{t,x;u^\varepsilon}-x|\geq \bar{\alpha}\}\wedge
(t+\alpha).$$  We observe that, for all $(s,y,z)\in [t,\tau]\times
{\mathbb{R}}\times {\mathbb{R}}^d$, from (\ref{F2}), (\ref{G1}) and the
definition of $\tau$
$$\begin{array}{llll}
F(s,X_s^{t,x;u^\varepsilon},y,z, u_s^{\varepsilon}) &\leq & C^*(\vert
y\vert+\vert
z\vert)+F(s,X_s^{t,x;u^\varepsilon},0,0,u_s^{\varepsilon})\\
& \leq & C^*(\vert y\vert+\vert z\vert)-\frac{1}{2}\theta;
       \end{array}$$
 $$\begin{array}{llll}
G(s,X_s^{t,x;u^\varepsilon},y) \leq C^*\vert
y\vert+G(s,X_s^{t,x;u^\varepsilon},0) \leq C^*\vert y\vert-\frac{1}{2}\theta.
       \end{array}$$
Consequently, applying the comparison result for GBSDEs (Lemma 2.7, or
Remark 1.5 in Pardoux and Zhang~\cite{PZ}) to GBSDEs (4.6) and (\ref{4.21}) we have that
\be Y_s^{1,u^\varepsilon}\leq Y^4_s,\, s\in[t,\tau],\ \mbox{P-a.s.,}\label{4.20}\ee where $Y^4$ is
defined by the following BSDE: \be\label{4.21} \left \{\begin{array}{lll}
      &-dY^4_s =(C^*(\vert Y^4_s\vert+\vert Z^4_s\vert)-\frac{1}{2}\theta)ds+
      (C^*|Y_s^4|-\frac{1}{2}\theta)dK_s^{t,x;u^\varepsilon}-Z^4_s dB_s, \\
     &Y^4_{\tau}=Y_\tau^{1;u^\varepsilon}.
     \end{array}\right.  \ee
On the other hand, we also have to introduce the following BSDE: \be\label{4.22}
\left \{\begin{array}{lll}
      &-dY^5_s =(C^*(\vert Y^5_s\vert+\vert Z^5_s\vert)-\frac{1}{2}\theta)ds+
      (C^*|Y_s^5|-\frac{1}{2}\theta)dK_s^{t,x;u^\varepsilon}-Z^5_s dB_s, \\
     &Y^5_{t+\alpha}=0.
     \end{array}\right.  \ee
Notice that $C^*\vert Y^2_s\vert-\frac{1}{2}\theta<0,$ therefore
$Y^5_s\leq Y^2_s,\ s\in [t,t+\alpha],\ \mbox{P-a.s.}$, from the comparison theorem-Lemma 2.4. From Lemma 2.3 we have
$$\begin{array}{rcl} |Y_t^4-Y_t^5|&\leq&C(E[|Y_{\tau}^{4}-Y_{\tau}^{5}|^2|{\cal{F}}_t])^{\frac{1}{2}}\leq C(E[|Y_{\tau}^{4}|^2|{\cal{F}}_t])^{\frac{1}{2}}
+C(E[|Y_{\tau}^{5}|^2|{\cal{F}}_t])^{\frac{1}{2}}\\
&\leq& C(E[|Y_{\tau}^{1;u^\varepsilon}|^2|{\cal{F}}_t])^{\frac{1}{2}}
+C(E[|Y_{\tau}^{2}|^2|{\cal{F}}_t])^{\frac{1}{2}}+C(E[|Y_{\tau}^{5}-Y_{\tau}^{2}|^2|{\cal{F}}_t])^{\frac{1}{2}}\\
&\leq&  C\alpha^{\frac{3}{2}}+C(E[|Y_{\tau}^{5}-Y_{\tau}^{2}|^2|{\cal{F}}_t])^{\frac{1}{2}}\ \mbox{(from the proof of (\ref{4.16}))},\\
\end{array}
$$ for any $\alpha\in (0,\bar{\alpha}].$\ From (\ref{5.17-1}) of Remark 5.3 in Appendix, similarly we also have
\be\label{4.23}
P\{\tau<t+\alpha\, |{\cal{F}}_t\}\leq P\{\sup_{s\in [t,
t+\alpha]}|{X}_s^{t,x;u^\varepsilon}-x|\geq\bar{\alpha}|{\cal{F}}_t\}\leq \frac{C}{\bar{\alpha}^8}\alpha^4.  \ee
On the other hand, from Lemma 2.3 (taking $\mu=1$)\be\label{4.24}\begin{array}{lll}
&E[|Y_{\tau}^{5}-Y_{\tau}^{2}|^2|{\cal{F}}_t]\leq CE[\int^{t+\alpha}_\tau e^{2 K_s^{t,x;u^\varepsilon}}  (C^*|Y_s^2|-\frac{1}{2}\theta)^2dK_s^{t,x;u^\varepsilon}|{\cal{F}}_t]\\
=&CE[\int^{t+\alpha}_\tau e^{2 K_s^{t,x;u^\varepsilon}}  \frac{\theta^2}{4}e^{2C^*(s-(t+\alpha))}dK_s^{t,x;u^\varepsilon}|{\cal{F}}_t]\\
\leq & C \frac{\theta^2}{4}E[I_{\{\tau<t+\alpha\}}(e^{2K_{t+\alpha}^{t,x;u^\varepsilon}}-e^{2K_\tau^{t,x;u^\varepsilon}})|{\cal{F}}_t]\\
\leq & C \frac{\theta^2}{4}E[I_{\{\tau<t+\alpha\}}e^{2K_{t+\alpha}^{t,x;u^\varepsilon}}(K_{t+\alpha}^{t,x;u^\varepsilon}-K_\tau^{t,x;u^\varepsilon})|{\cal{F}}_t]\\
\leq &C\frac{\theta^2}{4}(P[\tau<t+\alpha|{\cal{F}}_t])^{\frac{1}{4}}(E[e^{8K_{t+\alpha}^{t,x;u^\varepsilon}}|{\cal{F}}_t])^{\frac{1}{4}}
(E[|K_{t+\alpha}^{\tau,X_\tau^{t,x;u^\varepsilon};u^\varepsilon}|^2|{\cal{F}}_t])^{\frac{1}{2}}\\
\leq &
C\frac{\theta^2}{4}(P\{\tau<t+\alpha|{\cal{F}}_t\})^{\frac{1}{4}}(E[(t+\alpha-\tau)^2|{\cal{F}}_t])^{\frac{1}{4}}\ \hfill\mbox{(from Prop. 5.1 and 5.3.)}\\
\leq &C\theta^2(P\{\tau<t+\alpha|{\cal{F}}_t\})^{\frac{1}{4}}(\alpha^2P\{\tau<t+\alpha\, |{\cal{F}}_t\})^{\frac{1}{4}}\\
\leq &C \theta^2\alpha^{\frac{5}{2}}.\\
\end{array}
\ee

Therefore, \be\label{4.25} |Y_t^4-Y_t^5|\leq C\alpha^{\frac{3}{2}}+C\theta
\alpha^{\frac{5}{4}}.\ee

Now we obtain $$- \varepsilon\alpha\leq Y^{1,u^\varepsilon}_t\leq Y_t^4\leq
Y_t^5+|Y_t^4-Y_t^5|\leq Y_t^2+C\alpha^{\frac{3}{2}}+C\theta
\alpha^{\frac{5}{4}}, \mbox{P-a.s.}$$ i.e.,\ \ \ $- \varepsilon\alpha\leq Y^{1,u^\varepsilon}_t\leq - \frac{\theta }{2C^*}(1-e^{-
C^*\alpha})+C\alpha^{\frac{3}{2}}+C\theta \alpha^{\frac{5}{4}},\ \
\mbox{P-a.s.}$\ \ \  Therefore, $$- \varepsilon\leq - \frac{\theta }{2 C^*}\frac{1-e^{-
C^*\alpha}}{\alpha}+C\alpha^{\frac{1}{2}}+C\theta
\alpha^{\frac{1}{4}}, $$ and by taking the limit as
$\alpha\downarrow 0, \varepsilon\downarrow 0$ we get $0\leq -
\frac{\theta }{2}$ which contradicts our assumption that $\theta>0$.
Therefore, it must hold
$$\max\{\frac{\partial \varphi}{\partial t} (t,x) +  H(t, x,
W, D\varphi, D^2\varphi), \ \ \frac{\partial \varphi}{\partial n}
(t,x) + g(t, x, W)\} \geq 0.$$
\endpf

In an identical way, we can show
\begin{proposition}\label{pr4.2}Under the assumptions {\rm (H4.1)} and {\rm (H4.2)},
 the value function $W$ is a viscosity supersolution to {\rm (\ref{4.1})}.
\end{proposition} \noindent

\noindent{\it Proof}. Obviously, $W(T,x)=\Phi (x),\ x\in \bar{D}$. Suppose
that $\varphi \in C^3_{l,b} ([0,T] \times \bar{D })$ and $(t,x)\in
[0, T)\times \bar{D }$ is such that $W-\varphi$\ attains its minimum
at $(t,x)$. Without loss of generality, assume that $\varphi
(t,x)=W(t,x)$.

We first consider the case $x\in D$. We shall prove that
$$      \sup_{u\in U}F(t,x,0,0,u)\leq 0. $$
If this is not true, then there exists some $\theta>0$ such that
\be\label{4.26} F_0(t,x):=\sup_{u\in U}
F(t,x,0,0,u)\geq \theta>0.\ee Therefore, there exists a $u^{*}=u^{*}(t,x)\in U$ such that\ $ F(t,x,0,0,u^{*})\geq \frac{2\theta}{3}.$

Since $F_0$ is continuous at $(t,x)$, we can choose $\bar{\alpha}\in
(0, T-t]$ (for simplifying the notation, we still use $\bar{\alpha}$) such that
\begin{eqnarray}
&O_{\bar{\alpha}}(x):=\{y:|y-x|\leq \bar{\alpha}\}\subset D, \\
&F(s,y,0,0,u^{*})\geq \frac{1}{2}\theta \mbox{ \rm for all }(s,y)\in
[t, t+\bar{\alpha}]\times O_{\bar{\alpha}}(x).\label{F}
\end{eqnarray}
For any $\alpha\in (0, \bar{\alpha}]$, we still consider the BSDE (\ref{4.4}):
    \be \left \{\begin{array}{rcl}
      -dY^{1,u}_s &=& F(s,X^{t,x;u}_s,
      Y^{1,u}_s,Z^{1,u}_s,u_s)\, ds+G(s,X^{t,x;u}_s,
      Y^{1,u}_s)\, dK_s^{t,x;u}\\
      && -Z^{1,u}_s dB_s, \quad s\in [t,t+\alpha]; \\
     Y^{1,u}_{t+\alpha}&=&0,
     \end{array}\right.\ee
     where the pair of processes $(X^{t,x,u},K^{t,x,u})$ are given by
     $(\ref{3.1})$\ and $u(\cdot) \in {\mathcal{U}}_{t, t+\alpha}$. Therefore, Lemma 4.1 still holds for $x\in \bar{D}$.
On the other hand, from the DPP (Theorem 3.1), for
every $\alpha$,
     $$
     \varphi (t,x) =W(t,x) =\esssup_{u \in
{\mathcal{U}}_{t, t+\alpha}}G^{t,x;u}_{t,t+\alpha} [W(t+\alpha,
X^{t,x;u}_{t+\alpha})],
     $$
  and from $W\geq\varphi$\ and the monotonicity property of $G^{t,x;u}_{t,t+\delta}[\cdot]$\ (see Lemma 2.4)\ we have
     $$
     \esssup_{u \in
{\mathcal{U}}_{t, t+\alpha}} \{G^{t,x;u}_{t,t+\delta}
[\varphi(t+\alpha, X^{t,x;u}_{t+\alpha})] -\varphi (t,x)\}\leq 0,\
\mbox{{\it P}-a.s.}
     $$
 Thus, from Lemma \ref{l4.1}, we get $\label{4.29}
    \esssup_{u \in
{\mathcal{U}}_{t, t+\alpha}} Y^{1,u}_t \leq 0,\ \ \mbox{{\it P}-a.s.}, $\
which implies that \be\label{4.30}
Y^{1,u^{*}}_t\leq 0,\ \ \mbox{{\it P}-a.s.}
\ee
\br Similarly, the inequality (\ref{4.30}) holds true for $x\in \partial D$.
\er
For $u^{*}\in {\mathcal{U}}_{t, t+\alpha}$\ we define $\tau=\inf\{s\geq t: |X_s^{t,x;u^{*}}-x|\geq \bar{\alpha}\}\wedge
(t+\alpha).$\  Consequently, on $[t, \tau]$\ the process
$(K^{t,x;u^{*}})$\ is zero and, hence
$$Y_s^{1;u^{*}}=Y_\tau^{1;u^{*}}+\int_s^\tau F(r,X_r^{t,x;u^{*}},Y_r^{1;u^{*}},
Z_r^{1;u^{*}}, u^{*})dr-\int_s^\tau
Z_r^{1;u^{*}}dB_r. $$ We consider the following two BSDEs:
\be\label{4.33} \left \{\begin{array}{rcl}
      -d\widehat{Y}^2_s &=& (-C^*(\vert \widehat{Y}^2_s\vert+\vert \widehat{Z}^2_s\vert)+\frac{1}{2}\theta)\, ds-\widehat{Z}^2_s\, dB_s, \\
     \widehat{Y}^2_{t+\alpha}&=&0,
     \end{array}\right.
     \ee
whose unique solution is given by \be\label{4.34}
\widehat{Y}^2_s=\frac{\theta}{2C^*}(1-e^{C^*(s-(t+\alpha))}),\ \
\widehat{Z}^2_s=0, \quad s\in[t,t+\alpha],\ee
and \be \label{4.35}\left
\{\begin{array}{rcl}
      -d\widehat{Y}^3_s &=& (-C^*(\vert \widehat{Y}^3_s\vert+\vert \widehat{Z}^3_s\vert)+\frac{1}{2}\theta)\, ds-\widehat{Z}^3_s\, dB_s, \quad s\in[t,\tau];\\
     \widehat{Y}^3_{\tau}&=& Y_\tau^{1;u^*}.
     \end{array}\right.
     \ee
We have the following lemma.

\bl We have $ Y_t^{1,u^*}
 \geq \widehat{Y}_t^3$ and $|\widehat{Y}_t^2-\widehat{Y}_t^3|\leq C\alpha^{\frac{3}{2}},
 \mbox{P-a.s.}$\
 Here $C>0$\ is  independent of both the control $u$ and
$\alpha$.\el

\noindent{\it Proof}. (1) We observe from (\ref{F}) and the definition of $\tau$
that, for all $(s,y,z,u)\in [t,\tau]\times {\mathbb{R}}\times
{\mathbb{R}}^d\times U$,
$$
       \begin{array}{llll}
F(s,X_s^{t,x;u^*},y,z, u^*) &\geq & -C^*(\vert y\vert+\vert
z\vert)+F(s,X_s^{t,x;u^*},0,0,u^*)  \geq  -C^*(\vert y\vert+\vert z\vert)+\frac{1}{2}\theta.
       \end{array}$$
Consequently, from Lemma 2.2 in \cite{BL} we have that\ $Y_s^{1,u^*}\geq \widehat{Y}^3_s,\ s\in[t,\tau],$\
 where $\widehat{Y}^3$ is defined by BSDE~(\ref{4.35}).

(2) From the equation (4.31), Propositions 5.1 and 5.2
$$|Y_\tau^{1;u^*}|\leq C(t+\alpha-\tau)^{\frac{1}{2}}
+C(E[(K_{t+\alpha}^{t,x;u^*}-K_{\tau}^{t,x;u^*})^2|{\cal{F}}_\tau])^{\frac{1}{2}},$$ where $C$ is independent of controls. Then similar to the proof of estimate (\ref{4.14-2}), we have
\be E[|Y_{\tau}^{1;u^*}|^2|{\cal{F}}_t]\leq
C(E[(t+\alpha-\tau)^2|{\cal{F}}_t])^{\frac{1}{2}}.\ee
Similar to (\ref{4.14}), we still have
 \be\label{4.37}
P\{\tau<t+\alpha\, |{\cal{F}}_t\}\leq \frac{C}{\bar{\alpha}^8}\alpha^4.\ee
Therefore, \be E[|Y_{\tau}^{1;u^*}|^2|{\cal{F}}_t]\leq
C\alpha (P\{\tau<t+\alpha|{\cal{F}}_t\})^{\frac{1}{2}} \leq
\frac{C}{\bar{\alpha}^4}\alpha^3.\ee Furthermore, from Lemma 2.3 in \cite{BL},
\be\label{4.39}\begin{array}{rcl} |\widehat{Y}_t^2-\widehat{Y}_t^3|&\leq&
C(E[|\widehat{Y}_{\tau}^{2}-\widehat{Y}_{\tau}^{3}|^2|{\cal{F}}_t])^{\frac{1}{2}}\leq C(E[|\widehat{Y}_{\tau}^{2}|^2|{\cal{F}}_t])^{\frac{1}{2}}
+C(E[|\widehat{Y}_{\tau}^{3}|^2|{\cal{F}}_t])^{\frac{1}{2}}\\
&\leq& C
\frac{\theta}{2}(1-e^{-C^*\alpha})(P\{\tau<t+\alpha|{\cal{F}}_t\})^{\frac{1}{2}}
+C(E[|Y_{\tau}^{1;u^*}|^2|{\cal{F}}_t])^{\frac{1}{2}}\\
&\leq&
C\frac{\theta}{2}(1-e^{-C^*\alpha})\frac{1}{\bar{\alpha}^4}\alpha^2+\frac{C}{\bar{\alpha}^2}\alpha^{\frac{3}{2}}\leq C\alpha^{\frac{3}{2}},\\
\end{array}
\ee for any $\alpha\in (0,\bar{\alpha}].$
\bigskip

\noindent\textbf{Proof of Proposition \ref{pr4.2} (sequel)}.

By combining (\ref{4.30}) with Lemma 4.3 we then obtain

$$0\geq Y^{1,u^*}_t\geq \widehat{Y}_t^3\geq \widehat{Y}_t^2-|\widehat{Y}_t^2-\widehat{Y}_t^3|\geq \widehat{Y}_t^2-C\alpha^{\frac{3}{2}}, \mbox{P-a.s.}$$
i.e., \ \ $0\geq Y^{1,u^*}_t\geq  \frac{\theta }{2C^*}(1-e^{-
C^*\alpha})-C\alpha^{\frac{3}{2}},\ \ \mbox{P-a.s.}$\ \ Therefore,
$$0\geq \frac{\theta }{2 C^*}\frac{1-e^{-
C^*\alpha}}{\alpha}-C\alpha^{\frac{1}{2}}.$$ Letting $\alpha\to 0+$
,  we get $0\geq \frac{\theta }{2}$, which
contradicts our assumption that $\theta>0$. Therefore, we have\ $
     \sup_{u\in U} F(t,x,0,0,u) \leq 0, $\      which implies by the definition of $F$  that $$\frac{\partial \varphi}{\partial t} (t,x) +  H(t, x,
W, D\varphi, D^2\varphi) \leq 0,\ \  \mbox{if}\ \ x\in
    D.$$

We now consider the case $x\in \partial D$. We must prove that
$$\min\{\frac{\partial \varphi}{\partial t} (t,x) +  H(t, x,
\varphi, D\varphi, D^2\varphi), \ \ \frac{\partial \varphi}{\partial
n} (t,x) + g(t, x, \varphi)\} \leq 0$$ If this is not true,
then there exists some $\theta>0$ such that \be\label{4.40}
\sup_{u\in
U} F(t,x,0,0,u)\geq\theta> 0,\quad
G(t,x,0)\geq\theta>0,\ee therefore, there exists $u^*\in U$\ such that $ F(t,x,0,0,u^*)\geq \frac{2\theta}{3}.$

Choose $\bar{\alpha}\in (0, T-t]$ such that \begin{eqnarray}
&F(s,y,0,0,u^*)\geq\frac{1}{2}\theta,\label{4.41-1} \\
&G(s,y,0)\geq \frac{1}{2}\theta,\  \mbox{for all}\   t\leq
s\leq t+\bar{\alpha},\ |y-x|\leq \bar{\alpha}.\label{4.42-1}
\end{eqnarray}

Now we fix $\bar{\alpha}$, and we consider any $\alpha\in
(0, \bar{\alpha}].$\ Similarly, we still consider GBSDE (4.31) with $x\in \partial D$. For this $u^*\in {\mathcal{U}}_{t,
t+\alpha}$\ we still have (\ref{4.30}) and define $$\tau=\inf\{s\geq t: |X_s^{t,x;u^*}-x|\geq \bar{\alpha}\}\wedge
(t+\alpha).$$  We observe that, for all $(s,y,z)\in [t,\tau]\times
{\mathbb{R}}\times {\mathbb{R}}^d$, from (\ref{4.41-1}), (\ref{4.42-1}) and the definition of $\tau$
$$\begin{array}{llll}
F(s,X_s^{t,x;u^*},y,z, u_s^{\varepsilon}) &\geq & -C^*(\vert
y\vert+\vert
z\vert)+F(s,X_s^{t,x;u^*},0,0,u_s^{\varepsilon})\\
& \geq & -C^*(\vert y\vert+\vert z\vert)
+\frac{1}{2}\theta;
       \end{array}$$
 $$\begin{array}{llll}
G(s,X_s^{t,x;u^*},y)\geq -C^*\vert
y\vert+G(s,X_s^{t,x;u^*},0) \geq -C^*\vert y\vert+\frac{1}{2}\theta.
       \end{array}$$
Consequently, from the comparison result for GBSDEs (Lemma 2.7, or
Remark 1.5 in~\cite{PZ}) we have that $Y_s^{1,u^*}\geq \widehat{Y}^4_s,\, s\in[t,\tau],\ \mbox{P-a.s.,}$\ where $\widehat{Y}^4$ is
defined by the following BSDE: \be\label{4.41} \left \{\begin{array}{lll}
      &-d\widehat{Y}^4_s =(-C^*(\vert \widehat{Y}^4_s\vert+\vert \widehat{Z}^4_s\vert)+\frac{1}{2}\theta)ds+
      (-C^*|\widehat{Y}_s^4|+\frac{1}{2}\theta)dK_s^{t,x;u^*}-\widehat{Z}^4_s dB_s, \\
     &\widehat{Y}^4_{\tau}=Y_\tau^{1;u^*}.
     \end{array}\right.  \ee
On the other hand, we also have to introduce the following BSDE: \be\label{4.42}
\left \{\begin{array}{lll}
      &-d\widehat{Y}^5_s =(-C^*(\vert \widehat{Y}^5_s\vert+\vert \widehat{Z}^5_s\vert)+\frac{1}{2}\theta)ds+
      (-C^*|\widehat{Y}_s^5|+\frac{1}{2}\theta)dK_s^{t,x;u^*}-\widehat{Z}^5_s dB_s, \\
     &\widehat{Y}^5_{t+\alpha}=0.
     \end{array}\right.  \ee
Notice that $-C^*\vert \widehat{Y}^2_s\vert+\frac{1}{2}\theta>0,$\ therefore
$\widehat{Y}^5_s\geq \widehat{Y}^2_s,\ s\in [t,t+\alpha],\ \mbox{P-a.s.}$, from Lemma 2.4.

From Lemma 2.3 we have \be\begin{array}{rcl} |\widehat{Y}_t^4-\widehat{Y}_t^5|&\leq&
C(E[|\widehat{Y}_{\tau}^{4}-\widehat{Y}_{\tau}^{5}|^2|{\cal{F}}_t])^{\frac{1}{2}}\leq C(E[|\widehat{Y}_{\tau}^{4}|^2|{\cal{F}}_t])^{\frac{1}{2}}
+C(E[|\widehat{Y}_{\tau}^{5}|^2|{\cal{F}}_t])^{\frac{1}{2}}\\
&\leq& C(E[|Y_{\tau}^{1;u^*}|^2|{\cal{F}}_t])^{\frac{1}{2}}
+C(E[|\widehat{Y}_{\tau}^{2}|^2|{\cal{F}}_t])^{\frac{1}{2}}+C(E[|\widehat{Y}_{\tau}^{5}-\widehat{Y}_{\tau}^{2}|^2|{\cal{F}}_t])^{\frac{1}{2}}\\
&\leq&  C\alpha^{\frac{3}{2}}+C(E[|\widehat{Y}_{\tau}^{5}-\widehat{Y}_{\tau}^{2}|^2|{\cal{F}}_t])^{\frac{1}{2}}\hfill \mbox{(from the proof of (\ref{4.39}))},\\
\end{array}
\ee for any $\alpha\in (0,\bar{\alpha}].$

Similar to (4.25) and (4.26),\ $P\{\tau<t+\alpha\, |{\cal{F}}_t\}\leq\frac{C}{\bar{\alpha}^8}\alpha^4;$\ and
\be\begin{array}{lll}
&E[|\widehat{Y}_{\tau}^{5}-\widehat{Y}_{\tau}^{2}|^2|{\cal{F}}_t]\leq CE[\int^{t+\alpha}_\tau e^{2 K_s^{t,x;u^*}}  (C^*|\widehat{Y}_s^2|-\frac{1}{2}\theta)^2dK_s^{t,x;u^*}|{\cal{F}}_t]\\
=&CE[\int^{t+\alpha}_\tau e^{2 K_s^{t,x;u^*}} \frac{\theta^2}{4}e^{2C^*(s-(t+\alpha))}dK_s^{t,x;u^*}|{\cal{F}}_t]\leq C \theta^2\alpha^{\frac{5}{2}}.
\end{array}
\ee
Therefore, \be |\widehat{Y}_t^4-\widehat{Y}_t^5|\leq C\alpha^{\frac{3}{2}}+C\theta
\alpha^{\frac{5}{4}}.\ee

Now we obtain
$$0\geq Y^{1,u^*}_t\geq \widehat{Y}_t^4\geq
\widehat{Y}_t^5-|\widehat{Y}_t^4-\widehat{Y}_t^5|\geq \widehat{Y}_t^2-C\alpha^{\frac{3}{2}}-C\theta
\alpha^{\frac{5}{4}},\ \mbox{P-a.s.}$$ i.e.,\ \ $0\geq Y^{1,u^*}_t\geq \frac{\theta }{2C^*}(1-e^{-
C^*\alpha})-C\alpha^{\frac{3}{2}}-C\theta \alpha^{\frac{5}{4}},\ \mbox{P-a.s.}$\ \ Therefore,
$$0\geq  \frac{\theta }{2 C^*}\frac{1-e^{-
C^*\alpha}}{\alpha}-C\alpha^{\frac{1}{2}}-C\theta
\alpha^{\frac{1}{4}}\, ,$$ and by taking the limit as
$\alpha\downarrow 0, $ we get $0\geq \frac{\theta }{2}$ which contradicts our assumption that $\theta>0$.
Therefore, it must hold
$$\min\{\frac{\partial \varphi}{\partial t} (t,x) +  H(t, x,
W, D\varphi, D^2\varphi), \ \ \frac{\partial \varphi}{\partial n}
(t,x) + g(t, x, W)\} \leq 0.$$
\endpf

Therefore, we have

\bt Under Assumptions {\rm (H4.1)} and {\rm
(H4.2)}, the value function $W$ is the unique viscosity solution to {\rm
(\ref{4.1})}. \et

\br From Propositions~\ref{pr4.1} and \ref{pr4.2}, it remains to show the uniqueness assertion, which can be referred to Barles~\cite[Section 3]{B}, Bourgoing~\cite[Section 3]{BO}, and Crandall, Ishii, and
Lions~\cite[Section 7B]{CIL}.\er

\section{\large{Appendix}}

\subsection{\large{Forward-Backward SDES (FBSDEs)}}

 \hskip0.7cm In this section we give some necessary basic results on GBSDEs associated
 with forward reflected SDEs (for short: FSDEs). We consider measurable functions $b:[0,T]\times \Omega\times
{\mathbb{R}}^d\rightarrow {\mathbb{R}}^d \ $ and
         $\sigma:[0,T]\times \Omega\times {\mathbb{R}}^d\rightarrow {\mathbb{R}}^{d\times d}$
which are supposed to satisfy the following conditions:
$$\begin{array}{lll}
&{\rm{(i)}}\  b(\cdot,0)\ \mbox{and}\ \sigma(\cdot,0)\ \mbox{are} \ {\mathbb{F}}-\mbox{adapted processes; there exists some constant}\ C>0,\ \mbox{such that}\\
&\ \ \ |b(t,x)|+|\sigma(t,x)|\leq C(1+|x|),\ \mbox{a.s.,}\ \mbox{for all}\ 0\leq t\leq T,\ x\in {\mathbb{R}}^d;\\
&{\rm{(ii)}}\ b\ \mbox{and}\ \sigma\ \mbox{are Lipschitz in}\ x,\ \mbox{i.e., there is some constant}\ C>0\ \mbox{such that}\hfill{{\rm (H5.1)}}\\
&\ \ \ |b(t,x)-b(t,x')|+|\sigma(t,x)-\sigma(t,x')|\leq C| x-x'|,\ \mbox{a.s.,}\ \mbox{for all}\ 0\leq t \leq T,\ x,\ x'\in {\mathbb{R}}^d.\\
\end{array}$$

Under the assumption (H5.1), it follows from the results in Lions
and Sznitman~\cite{LS} that for each initial condition $(t,\zeta)\in[0,T]\times L^2(\Omega,{\cal{F}}_t,P;\bar{D})$\ there
exists a unique pair of progressively measurable continuous processes $\{(X^{t,\zeta},K^{t,\zeta})\}$, with values in
$\bar{D}\times {\mathbb{R}}_+$, such that
  \be
  \left\{
  \begin{array}{rcl}
  X_s^{t,\zeta}&=&\zeta+\int_t^sb(r,X_r^{t,\zeta})dr+\int_t^s\sigma(r,X_r^{t,\zeta})dB_r+\int_t^s\nabla\phi(X_r^{t,\zeta})dK_r^{t,\zeta},\ s\in[t,T],\\
  K_s^{t,\zeta}&=&\int_t^sI_{\{X_r^{t,\zeta}\in\partial D\}}dK_r^{t,\zeta},\ K^{t,\zeta}\ \mbox{is increasing}.
  \end{array}
  \right.
  \ee
\bp For each $T\geq 0$, there exists a constant $C_T$\ such that, for
all $\zeta,\ \zeta'\in L^2(\Omega,{\cal{F}}_t,P;\bar{D})$, \be
E(\sup_{t\leq s\leq
T}|X_s^{t,\zeta}-X_s^{t,\zeta'}|^4|{\cal{F}}_t)\leq
C_T|\zeta-\zeta'|^4,\ee and \be E(\sup_{t\leq s\leq
T}|K_s^{t,\zeta}-K_s^{t,\zeta'}|^4|{\cal{F}}_t)\leq
C_T|\zeta-\zeta'|^4. \ee Moreover, for each $\mu>0, s\in [t,T]$,
there exists $C(\mu, s)$ such that for all $\zeta\in
L^2(\Omega,{\cal{F}}_t,P;\bar{D}),$ \be E(e^{\mu
K_s^{t,\zeta}}|{\cal{F}}_t)\leq C(\mu, s). \ee \ep

The proof is similar to that of Propositions 3.1 and  3.2 in Pardoux
and Zhang~\cite{PZ}.

We assume that the three functions $f,\ g$ and $\Phi$\ satisfy the following conditions:
$$\begin{array}{lll}
&{\rm{(i)}}\ \Phi: \Omega\times {\mathbb{R}}^d\rightarrow {\mathbb{R}} \ \mbox{is an}\ {\cal{F}}_T\otimes{\cal{B}}({\mathbb{R}}^d)
             \mbox{-measurable random variable and}\ f: [0,T]\times \Omega\times {\mathbb{R}}^d\\
& \ \ \ \times {\mathbb{R}}\times {\mathbb{R}}^d \rightarrow {\mathbb{R}}\ \mbox{is a measurable process}\ \mbox{such that}\ f(\cdot,x,y,z)\ \mbox{is}\ {\mathbb{F}} \mbox{-adapted, for all }(x, y, z)\\
& \ \ \ \in{\mathbb{R}}^d\times {\mathbb{R}}\times {\mathbb{R}}^d;\ g: [0,T]\times {\mathbb{R}}^d\times {\mathbb{R}} \rightarrow {\mathbb{R}}\ \mbox{is a measurable function}\ \mbox{such that}\\
& \ \ \  g(\cdot)\in C^{1,2,2}([0,T]\times {\mathbb{R}}^d\times {\mathbb{R}});\\
&{\rm{(ii)}}\ \mbox{There exists a constant}\ C>0\ \mbox{such that}\ \mbox{for all}\ t\in[0, T],\ x, x', z, z'\in {\mathbb{R}}^d,\ y, y'\in {\mathbb{R}}, \\
& \ \ \ | f(t,x,y,z)-f(t,x',y',z')| +|g(t,x,y)-g(t,x',y')| +| \Phi(x)-\Phi(x')|\\
& \ \ \ \leq C(|x-x'|+ |y-y'|+|z-z'|),\ \mbox{ a.s.};\hfill{{\rm (H5.2)}}\\
&{\rm{(iii)}}\ f\ \mbox{and}\ \Phi \ \mbox{satisfy a linear growth condition, i.e., there exists some}\ C>0\ \ \mbox{such that},\\
& \ \ \ \mbox{for all}\ x\in {\mathbb{R}}^d, |f(t,x,0,0)|+ |\Phi(x)| \leq C(1+|x|),\ \mbox{a.s.a.e.}
    \end{array}$$

 Under the above assumptions the coefficients $f(s,X_s^{t,\zeta},y,z)$\ and $g(s,X_s^{t,\zeta},y)$\ satisfy (H2.1) and $\xi=\Phi(X_T^{t,\zeta})$ $\in
 L^2(\Omega,{\cal{F}}_T,P)$. Therefore, the following GBSDE possesses a unique solution:
\be\label{5.5} \left\{
\begin{array}{rcl}
-dY_s^{t,\zeta}&=&f(s,X_s^{t,\zeta},Y_s^{t,\zeta},Z_s^{t,\zeta})ds+g(s,X_s^{t,\zeta},Y_s^{t,\zeta})dK_s^{t,\zeta}-Z_s^{t,\zeta}dB_s,\ s\in [t, T],\\
Y_T^{t,\zeta}&=&\Phi(X_T^{t,\zeta}).\\
\end{array}
\right. \ee

\bp  Let assumptions (H5.1) and (H5.2) hold. Then, for any $0\leq t\leq T$ and $\zeta,\ \zeta'\in L^2(\Omega,{\cal{F}}_t,P;\bar{D})$,  $$
 \begin{array}{lll}
\mbox{\rm(i)}& E[\sup_{t\leq s\leq
T}|Y_s^{t,\zeta}|^2+\int_t^T|Z_s^{t,\zeta}|^2ds|{\cal{F}}_t]\leq
C(1+|\zeta|^2),\  a.s.; \mbox{ and in particular,}\\
&\ |Y_t^{t,\zeta}|\leq C(1+|\zeta|),\  a.s.; \hskip3cm\\
\mbox{\rm(ii)}&|Y_t^{t,\zeta}-Y_t^{t,\zeta'}|\leq C|\zeta-\zeta'|+C|\zeta-\zeta'|^{\frac{1}{2}},\  a.s.,\hskip3cm \\
\end{array}
$$ where the constant $C>0$\ depends only on the Lipschitz and the
growth constants of $b$,\ $\sigma$, $f$, $g$\ and $\Phi$. \ep

\br Since $D$ is bounded, we have \be E[\sup_{t\leq s\leq T}|Y_s^{t,\zeta}|^2+\int_t^T|Z_s^{t,\zeta}|^2ds|{\cal{F}}_t]\leq C,\
\mbox{a.s.},\ee where $C$ is independent of $\zeta$.  \er

\noindent{\it Proof.} From Lemma 2.2 and Proposition 5.1, we have assertion (i). Now we prove assertion (ii). First notice that from (i) we have $|Y^{t,\zeta}_t|\leq C (1+|\zeta|),$\ a.s., therefore we can get from the uniqueness of the solution of equations (5.1) and (5.5) that
\be |Y^{t,\zeta}_s|=|Y^{s,X_s^{t,\zeta}}_s|\leq C(1+|X_s^{t,\zeta}|)\leq C,\ \mbox{a.s.},\ee
since D is bounded. From Burkholder-Davis-Gundy inequality and (5.5), as well as from the boundedness of the processes $X^{t,\zeta},\ Y^{t,\zeta}$,
$$
 \begin{array}{lll}&E[(\int_s^T|Z^{t,\zeta}_r|^2dr)^2|{\cal F}_t] \le CE[\sup_{r\in[s,T]}|\int_s^rZ^{t,\zeta}_v dB_v|^4|{\cal F}_t]\\
&\le C+C_0(T-s)^2E[(\int_s^T|Z^{t,\zeta}_r|^2dr)^2|{\cal F}_t] +
CE[(K^{t,\zeta}_T)^4|{\cal F}_t]\\
&\le C+ C_0(T-s)^2E[(\int_s^T|Z^{t,\zeta}_r|^2dr)^2|{\cal F}_t].
\end{array}
$$
Consequently, for $T-s\le (\frac{1}{2C_0})^{1/2}$, $E[(\int_s^T|Z^{t,\zeta}_r|^2dr)^2|{\cal F}_t]\le C$. This argument allows to choose a partition $t=t_0<t_1<...<t_N=T$ of the interval $[t,T]$\ such that $E[(\int_{t_{i-1}}^{t_i}|Z^{t,\zeta}_r|^2dr)^2|{\cal F}_t]\le C$, $1\le i\le N.$\ Therefore, we have
\be E[(\int_t^T|Z^{t,\zeta}_r|^2dr)^2|{\cal F}_t]\le C.\ee
For any $\lambda>0$, applying
It\^{o}'s formula to $e^{\lambda K_s^{t,\zeta'}}|Y^{t,\zeta}_s-Y^{t,\zeta'}_s|^2$, we have \be
\begin{array}{lll}
&|Y^{t,\zeta}_s-Y^{t,\zeta'}_s|^2+\lambda\int_s^Te^{\lambda
K_r^{t,\zeta'}}|Y^{t,\zeta}_r-Y^{t,\zeta'}_r|^2dK_r^{t,\zeta'}+\int_s^Te^{\lambda
K_r^{t,\zeta'}}|Z^{t,\zeta}_r-Z^{t,\zeta'}_r|^2dr\\
=&e^{\lambda
K_T^{t,\zeta'}}|Y^{t,\zeta}_T-Y^{t,\zeta'}_T|^2-2\int_s^Te^{\lambda K_r^{t,\zeta'}}(Y^{t,\zeta}_r-Y^{t,\zeta'}_r)<Z^{t,\zeta}_r-Z^{t,\zeta'}_r,dB_r>\\
&+2\int_s^Te^{\lambda K_r^{t,\zeta'}}(Y^{t,\zeta}_r-Y^{t,\zeta'}_r)(
f(r,X^{t,\zeta}_r,Y^{t,\zeta}_r,Z^{t,\zeta}_r)-f(r,X^{t,\zeta'}_r,Y^{t,\zeta'}_r,Z^{t,\zeta'}_r))dr\\
&+2\int_s^Te^{\lambda K_r^{t,\zeta'}}(Y^{t,\zeta}_r-Y^{t,\zeta'}_r)
(g(r,X^{t,\zeta}_r,Y^{t,\zeta}_r)-g(r,X^{t,\zeta'}_r,Y^{t,\zeta'}_r))dK_r^{t,\zeta'}\\
&+2\int_s^Te^{\lambda K_r^{t,\zeta'}}(Y^{t,\zeta}_r-Y^{t,\zeta'}_r)
g(r,X^{t,\zeta}_r,Y^{t,\zeta}_r)d(K_r^{t,\zeta}-K_r^{t,\zeta'}).\\
\end{array}
\ee Then from (H5.1), (H5.2), (5.4), (5.7) and (5.8), taking a suitable $\lambda>0$, we  get \be
\begin{array}{lll}
 & |Y^{t,\zeta}_s-Y^{t,\zeta'}_s|^2 \leq C|\zeta-\zeta'|^2+CE[\int_s^Te^{\lambda
K_r^{t,\zeta'}}|Y_r^{t,\zeta}-Y_r^{t,\zeta'}|^2 dr|{\cal{F}}_s]\\
&\ \ \ +CE[\int_s^Te^{\lambda
K_r^{t,\zeta'}}|X_r^{t,\zeta}-X_r^{t,\zeta'}|^2 dr|{\cal{F}}_s]
+CE[\int_s^Te^{\lambda
K_r^{t,\zeta'}}|X_r^{t,\zeta}-X_r^{t,\zeta'}|^2
dK_r^{t,\zeta'}|{\cal{F}}_s]\\
&\ \ \ +E[2\int_s^Te^{\lambda
K_r^{t,\zeta'}}(Y^{t,\zeta}_r-Y^{t,\zeta'}_r)
g(r,X^{t,\zeta}_r,Y^{t,\zeta}_r)d(K_r^{t,\zeta}-K_r^{t,\zeta'})|{\cal{F}}_s].
\end{array}
\ee Furthermore, from Proposition 5.1, we have \be
\begin{array}{lll}
&|Y^{t,\zeta}_s-Y^{t,\zeta'}_s|^2 \leq C|\zeta-\zeta'|^2+CE[\int_s^Te^{\lambda
K_r^{t,\zeta'}}|Y_r^{t,\zeta}-Y_r^{t,\zeta'}|^2 dr|{\cal{F}}_s]\\
&\ \ \ +E[2\int_s^Te^{\lambda
K_r^{t,\zeta'}}(Y^{t,\zeta}_r-Y^{t,\zeta'}_r)
g(r,X^{t,\zeta}_r,Y^{t,\zeta}_r)d(K_r^{t,\zeta}-K_r^{t,\zeta'})|{\cal{F}}_s].
\end{array}
\ee
 On the other hand, applying It\^{o}'s formula to $e^{\lambda
K_s^{t,\zeta'}}(Y^{t,\zeta}_s-Y^{t,\zeta'}_s) g(s,X^{t,\zeta}_s,Y^{t,\zeta}_s)(K_s^{t,\zeta}-K_s^{t,\zeta'})$, we have \be
\begin{array}{lll}
 &E[\int_s^Te^{\lambda
K_r^{t,\zeta'}}(Y^{t,\zeta}_r-Y^{t,\zeta'}_r)
g(r,X^{t,\zeta}_r,Y^{t,\zeta}_r)d(K_r^{t,\zeta}-K_r^{t,\zeta'})|{\cal{F}}_s]\\
= & E[e^{\lambda K_T^{t,\zeta'}}(Y^{t,\zeta}_T-Y^{t,\zeta'}_T)
g(T,X^{t,\zeta}_T,Y^{t,\zeta}_T)(K_T^{t,\zeta}-K_T^{t,\zeta'})|{\cal{F}}_s] +E[\int_s^Tf_1(r)(K_r^{t,\zeta'}-K_r^{t,\zeta})dr|{\cal{F}}_s]\\
&+E[\int_s^Tf_2(r)(K_r^{t,\zeta'}-K_r^{t,\zeta})dK_r^{t,\zeta'}|{\cal{F}}_s] +E[\int_s^Tf_3(r)(K_r^{t,\zeta'}-K_r^{t,\zeta})dK_r^{t,\zeta}|{\cal{F}}_s],
\end{array}
\ee where
$$
\begin{array}{lll} f_1(s)=&-e^{\lambda
K_s^{t,\zeta'}}(f(s,X^{t,\zeta}_s,Y^{t,\zeta}_s,Z^{t,\zeta}_s)-f(s,X^{t,\zeta'}_s,Y^{t,\zeta'}_s,Z^{t,\zeta'}_s))
g(s,X^{t,\zeta}_s,Y^{t,\zeta}_s)\\
&+e^{\lambda
K_s^{t,\zeta'}}(Y_s^{t,\zeta}-Y_s^{t,\zeta'})\{\frac{\partial}{\partial
s}g(s,X_s^{t,\zeta},Y_s^{t,\zeta})+\nabla_xg(s,X_s^{t,\zeta},Y_s^{t,\zeta})b(s,X_s^{t,\zeta})\\
&-\nabla_yg(s,X_s^{t,\zeta},Y_s^{t,\zeta})f(s,X_s^{t,\zeta},Y_s^{t,\zeta},Z_s^{t,\zeta})
+\frac{1}{2}tr(D^2_xg(s,X_s^{t,\zeta},Y_s^{t,\zeta})\sigma\sigma^T(s,X_s^{t,\zeta}))\\
&+\frac{1}{2}D^2_yg(s,X_s^{t,\zeta},Y_s^{t,\zeta})|Z_s^{t,\zeta}|^2
+\frac{1}{2}tr<D_{xy}g(s,X_s^{t,\zeta},Y_s^{t,\zeta})
\sigma(s,X_s^{t,\zeta}),Z_s^{t,\zeta}>\}\\
&+e^{\lambda
K_s^{t,\zeta'}}(Z_s^{t,\zeta}-Z_s^{t,\zeta'})\{\nabla_xg(s,X_s^{t,\zeta},Y_s^{t,\zeta})\sigma(s,X_s^{t,\zeta})
+\nabla_yg(s,X_s^{t,\zeta},Y_s^{t,\zeta})Z^{t,\zeta}_s\};
\end{array}
$$
$f_2(s)=\{\lambda e^{\lambda
K_s^{t,\zeta'}}(Y_s^{t,\zeta}-Y_s^{t,\zeta'})+e^{\lambda
K_s^{t,\zeta'}}g(s,X_s^{t,\zeta'},Y_s^{t,\zeta'})\}g(s,X_s^{t,\zeta},Y_s^{t,\zeta});$
$$\begin{array}{lll} f_3(s)= &e^{\lambda
K_s^{t,\zeta'}}(Y_s^{t,\zeta}-Y_s^{t,\zeta'})
\{\nabla_xg(s,X_s^{t,\zeta},Y_s^{t,\zeta})\nabla
\phi(s,X_s^{t,\zeta})-
\nabla_yg(s,X_s^{t,\zeta},Y_s^{t,\zeta})g(s,X_s^{t,\zeta},Y_s^{t,\zeta})\}\\
&-e^{\lambda K_s^{t,\zeta'}}|g(s,X_s^{t,\zeta},Y_s^{t,\zeta})|^2.\\
\end{array}
$$
From assertion (i), Propositions 5.1, (5.7 )and (5.8),  we  have
$$E[\int_s^Te^{\lambda
K_r^{t,\zeta'}}(Y^{t,\zeta}_r-Y^{t,\zeta'}_r) g(r,X^{t,\zeta}_r,Y^{t,\zeta}_r)d(K_r^{t,\zeta}-K_r^{t,\zeta'})|{\cal{F}}_s]\leq
C|\zeta-\zeta'|^2+C|\zeta-\zeta'|.$$
Furthermore, from (5.11) and (5.4) we have $$\begin{array}{lll}
&E[|Y^{t,\zeta}_s-Y^{t,\zeta'}_s|^4|{\cal F}_t] \leq C|\zeta-\zeta'|^4+C|\zeta-\zeta'|^2+ CE[(\int_s^Te^{\lambda K_r^{t,\zeta'}}|Y^{t,\zeta}_r-Y^{t,\zeta'}_r|^2dr)^2|{\cal F}_t]\\
&\leq C|\zeta-\zeta'|^4+C|\zeta-\zeta'|^2+ C E[e^{2\lambda K_T^{t,\zeta'}}|{\cal F}_t]E[\int_s^T|Y^{t,\zeta}_r-Y^{t,\zeta'}_r|^4dr|{\cal F}_t]\\
&\leq C|\zeta-\zeta'|^4+C|\zeta-\zeta'|^2+ CE[\int_s^T|Y^{t,\zeta}_r-Y^{t,\zeta'}_r|^4dr|{\cal F}_t],\ \  s\in[t,T],
\end{array}
$$
then from Gronwall's Lemma, we get $E[|Y^{t,\zeta}_s-Y^{t,\zeta'}_s|^4|{\cal F}_t]\leq C|\zeta-\zeta'|^4+C|\zeta-\zeta'|^2,$ $ \mbox{a.s},\ s\in[t,T]$, which means (ii) for $s=t$.\endpf

\br If $g$\ is a bounded random variable, assertion {\rm(ii)} of (5.6)still holds. Indeed, from Lemma 2.3 in  \cite{BL} and Proposition 5.1,  we get
$$\begin{array}{rcl}
|Y_t^{t,\zeta}-Y_t^{t,\zeta'}|^2 &\leq &
CE[|\zeta-\zeta'+g(\omega)(K_T^{t,\zeta}-K_T^{t,\zeta'})|^2|{\cal{F}}_t]\\
& &+CE[\int_t^T|f(s,X_s^{t,\zeta},Y_s^{t,\zeta},Z_s^{t,\zeta})
-f(s,X_s^{t,\zeta'},Y_s^{t,\zeta},Z_s^{t,\zeta})|^2ds|{\cal{F}}_t]\\
&\leq & C|\zeta-\zeta'|^2,\ \mbox{a.s.} \end{array} $$\er

\bp Let assumptions (H5.1) and (H5.2) hold. Then, for any $0\leq \alpha\leq T-t$ and the associated initial conditions
 $\zeta \in L^2(\Omega,{\cal{F}}_t,P;\bar{D})$, we
have the following estimates:\be
E[|K_{t+\alpha}^{t,\zeta}|^2|{\cal{F}}_t]\leq
C\alpha,\ \mbox{a.s.}, \ee where the constant $C>0$\ depends only on
the Lipschitz and the growth constants of $b$,\ $\sigma$, $f$, $g$\
and $\Phi$. \ep

\noindent{\it Proof.} For $\zeta'\in L^2(\Omega,{\cal{F}}_t,P;\bar{D})$, from It\^{o}'s formula we have
\be
\begin{array}{rcl}
|X_s^{t,\zeta}-\zeta'|^2&=&|\zeta-\zeta'|^2+2\int_t^s(X_r^{t,\zeta}-\zeta')b(r,X_r^{t,\zeta})dr
 +2\int_t^s(X_r^{t,\zeta}-\zeta')\sigma(r,X_r^{t,\zeta})dB_r\\
 & &+\int_t^s|\sigma(r,X_r^{t,\zeta})|^2dr+2\int_t^s(X_r^{t,\zeta}-\zeta')\nabla\phi(X_r^{t,\zeta})dK^{t,\zeta}_r,\
 s\in[t, T].\end{array}
\ee Since $D\subset {\mathbb{R}}^d$\ is convex, we have \be\int_t^s(X_r^{t,\zeta}-\zeta')\nabla\phi(X_r^{t,\zeta})dK^{t,\zeta}_r\leq 0.\ee
Therefore, we have\ $ E[\sup_{s\in[t,t+\alpha]}|X_s^{t,\zeta}-\zeta'|^2|{\cal{F}}_t]\leq C(|\zeta-\zeta'|^2+\alpha). $\  Recall that $D$ is an
open connected bounded convex subset. In particular, we have, \be\label{5.15}
E[\sup_{s\in[t,t+\alpha]}|X_s^{t,\zeta}-\zeta|^2|{\cal{F}}_t]\leq C \alpha. \ee Because $\phi\in
 C_b^2({\mathbb{R}}^d)$\ we have
$$
\begin{array}{rcl}
\phi(X_s^{t,\zeta})&=&\phi(\zeta)+\int_t^s\nabla\phi(X_r^{t,\zeta})b(r,X_r^{t,\zeta})dr
 +\int_t^s\nabla\phi(X_r^{t,\zeta})\sigma(r,X_r^{t,\zeta})dB_r\\
 & &+\frac{1}{2}\int_t^str(D^2\phi\sigma(r,X_r^{t,\zeta})\sigma^T(r,X_r^{t,\zeta}))dr
 +\int_t^s|\nabla\phi(X_r^{t,\zeta})|^2dK^{t,\zeta}_r,\
 s\in[t, T].\end{array}
$$
Therefore, we  get
$$K^{t,\zeta}_s\leq |\phi(X_s^{t,\zeta})-\phi(\zeta)|+C\int_t^s(1+|X_r^{t,\zeta}|^2)dr
+|\int_t^s\nabla\phi(X_r^{t,\zeta})\sigma(r,X_r^{t,\zeta})dB_r|,$$ and furthermore, from Burkholder-Davis-Gundy inequality, we have
$$
E[|K^{t,\zeta}_{t+\alpha}|^2|{\cal{F}}_t]\leq CE[\sup_{s\in[t,t+\alpha]}|X_s^{t,\zeta}-\zeta|^2|{\cal{F}}_t]+C\alpha.
$$
In view of (\ref{5.15}), the proof is complete.\endpf

\br In view of (5.13) and (5.14),  using Burkholder-Davis-Gundy inequality, we  have \be\label{5.17-1}
E[\sup_{s\in[t,t+\alpha]}|X_s^{t,\zeta}-\zeta|^8|{\cal{F}}_t]\leq C \alpha^4. \ee \er

 Let us now define the random field:
\be u(t,x)=Y_s^{t,x}|_{s=t},\ (t, x)\in [0, T]\times\bar{D}, \ee
where $Y^{t,x}$ is the solution of GBSDE (\ref{5.5}) with $x \in \bar{D}$\
at the place of $\zeta\in
L^2(\Omega,{\cal{F}}_t,P;\bar{D}).$\\

Proposition 5.2 yields that, for all $t \in [0, T] $, P-a.s.,
 \be
\begin{array}{ll}
\mbox{(i)}&| u(t,x)-u(t,y)| \leq C|x-y|+C|x-y|^{\frac{1}{2}},\ \mbox{for all}\ x, y\in \bar{D};\\
\mbox{(ii)}&| u(t,x)|\leq C(1+|x|),\ \mbox{for all}\ x\in \bar{D}.\\
\end{array}
\ee
 \begin{theorem}\label{th6.1} Under the assumptions {\rm (H3.1)} and {\rm (H3.2)}, for any $t\in [0, T]$\ and $\zeta\in
L^2(\Omega,{\cal{F}}_t,P;\bar{D}),$\ we have \be\label{6.7}
u(t,\zeta)=Y_t^{t,\zeta},\ \mbox{{\it P}-a.s.} \ee \end{theorem}

The proof of Theorem \ref{th6.1} is similar to that of Theorem 3.1 in Peng~\cite{Pe1} or Theorem A.2 in Buckdahn and Li~\cite{BL}. Therefore
it is omitted here.

\subsection{Proofs of Proposition 3.1 and Theorem 3.1}
\noindent {\it Proof of Proposition 3.1}. Let $H$ be the Cameron--Martin space
 of all absolutely continuous elements $h\in \Omega$\ whose derivative  $\dot{h}$\ is in $L^2([0, T],{\mathbb{R}}^d).$

For any $h \in H$, we define $\tau_h\omega:=\omega+h,\
\omega\in \Omega$.  Obviously, $\tau_h: \Omega\rightarrow\Omega$ is
a bijection with the inverse $\tau_h^{-1}$.  The law is given by
$$P\circ[\tau_h^{-1}]=\exp\{\int^T_0\dot{h}_sdB_s-\frac{1}{2}\int^T_0|\dot{h}_s|^2ds\}P.$$
Fix any $(t, x)\in [0, T]\times \bar{D}$, and define $H_t: =\{h\in
H|h(\cdot)=h(\cdot\wedge t)\}.$\ The rest of the proof is divided
into the following three steps:

Step 1. For any $u\in {\mathcal{U}}_{t,T} \ \mbox{and}\ h \in H_t,\
J(t, x; u)$ $(\tau_h)= J(t, x;
u(\tau_h)),\ \mbox{{\it P}-a.s.}$ 

 Indeed, the $\tau_h$-shifted  reflected SDE~(\ref{3.1}) (with
 $\zeta=x$)  is the same reflected  SDE~(\ref{3.1}) with $u$ being substituted into
  the  $\tau_h$-shifted control
 process $u(\tau_h)$. From the uniqueness of the solution of the
 reflected SDE~(\ref{3.1}), we get $X_s^{t,x; u}(\tau_h)=X_s^{t,x;
 u(\tau_h)}$ and
$K_s^{t,x; u}(\tau_h)=K_s^{t,x; u(\tau_h)}$ for  $s\in [t, T]$ {\it
P}-a.s. Furthermore, by a similar shift argument and the associated
Girsanov transformation, we get from the uniqueness of the solution
of GBSDE~(\ref{3.5}) that
$$Y_s^{t,x; u}(\tau_h)=Y_s^{t,x; u(\tau_h)}\ \mbox{for
any}\ s\in [t, T],\ \mbox{{\it P}-a.s.,}$$
$$Z_s^{t,x; u}(\tau_h)=Z_s^{t,x; u(\tau_h)},\  \mbox{dsd{\it P}-a.e. on}\ [t, T]\times\Omega.$$
It means $$J(t, x; u)(\tau_h)= J(t, x; u(\tau_h)),\ \mbox{{\it
P}-a.s.}$$

Step 2. For all $h\in H_t$\ we have
$$\{\esssup_{u \in {\mathcal{U}}_{t,T}}J(t,x;
u)\}(\tau_h)=\esssup_{u \in {\mathcal{U}}_{t,T}}\{J(t,x;
u)(\tau_h)\},\ \mbox{{\it P}-a.s.}
$$

Indeed, define
 $$W(t,x):=\esssup_{u \in {\mathcal{U}}_{t,T}}J(t,x;
u).$$ we have  $W(t,x)\geq J(t,x; u),$\ and thus $W(t,x)(\tau_h)\geq
J(t,x; u)(\tau_h), \ \mbox{{\it P}-a.s.}, {\rm for~all}\ u \in
{\mathcal{U}}_{t,T}.$\ On the other hand, for any random variable
$\zeta$\ satisfying $\zeta\geq$ $ J(t,x; u)(\tau_h),$\ and hence
also $\zeta(\tau_{-h})\geq J(t,x; u), \ \mbox{{\it P}-a.s.,\ for}
\mbox{ all}\ u\in {\mathcal{U}}_{t,T},$\ we have\
$\zeta(\tau_{-h})\geq W(t,x), \ $ $ \mbox{{\it P}-a.s.,}$ i.e.,
$\zeta\geq W(t,x)(\tau_{h}),\mbox{{\it P}-a.s.}$\ Consequently,
$$W(t,x)(\tau_{h})=\esssup_{u \in
{\mathcal{U}}_{t,T}}\{J(t,x; u)(\tau_h)\},\ P\mbox{-a.s.}$$

Step 3. $W(t,x)$\ is invariant with respect
 to the shift $\tau_h$, i.e.,
  $$W(t,x)(\tau_{h})=W(t,x), \ \mbox{{\it P}-a.s., for any}\ h\in H. $$

Indeed, from Step 1 to Step 2, we have, for any $h\in H_t,$
 $$
   \begin{array}{rcl}
   W(t,x)(\tau_{h}) & = & \esssup_{u \in
{\mathcal{U}}_{t,T}}\{J(t,x; u)(\tau_h)\} = \esssup_{u \in
{\mathcal{U}}_{t,T}}J(t,x; u(\tau_h))\\
& = &  \esssup_{u \in
{\mathcal{U}}_{t,T}}J(t,x; u) = W(t,x),\ \mbox{{\it P}-a.s.,}
   \end{array}
$$
where we have used
$\{u(\tau_h)|u(\cdot)\in{\mathcal{U}}_{t,T}\}={\mathcal{U}}_{t,T} $\
so as to obtain the 3rd equality. Therefore, $W(t,x)(\tau_{h})=
W(t,x), \mbox{ {\it P}-a.s.}$ for any $h\in H_t$. Since $W(t,x)$ is
${\mathcal{F}}_{t}$-measurable, it holds for all $ h\in H.$ Indeed,
since $\Omega= C_0([0, T];{\mathbb{R}}^d)$, by the definition of the
filtration, the ${\cal F}_t$-measurable random variable $W(t,
x)(\omega),\ \omega\in \Omega,$\ only depends on the restriction of
$\omega$ to the time interval $[0, t]$.

 The result of Step 3,
combined with the following Lemma \ref{l3.1} (refer to Buckdahn and Li~\cite[Lemma 3.4]{BL}) completes the
proof.\endpf

\begin{lemma}\label{l3.1} Let $\zeta$\ be a random variable defined on
the Wiener space $(\Omega, {\mathcal{F}}_T, P)$ such that
 $\zeta(\tau_{h})=\zeta$  {\it P}-a.s. for any $h\in H.$ Then
$\zeta=E\zeta$ {\it P}-a.s.\end{lemma}

\bigskip

\noindent{\it Proof of Theorem 3.1}. To simplify our exposition, define
$$I_\delta(t, x, u):=G^{t,x;u}_{t,t+\delta} [W(t+\delta,
X^{t,x;u}_{t+\delta})] $$ and
$$W_\delta(t,x) :=\esssup_{u \in {\mathcal{U}}_{t,
t+\delta}}I_\delta(t, x, u)=\esssup_{u \in {\mathcal{U}}_{t,
t+\delta}}G^{t,x;u}_{t,t+\delta} [W(t+\delta,
X^{t,x;u}_{t+\delta})].$$
The proof of Theorem~\ref{th3.1} is reduced to the following three
lemmas. Similar to the proof of
Proposition~\ref{p3.1}, we first have

\begin{lemma}\label{l3.3} $W_\delta(t,x)$ is deterministic for
any $0\leq t<t+\delta \leq T,\ x\in \bar{D}$.\end{lemma}

\begin{lemma}\label{l3.4}$W_\delta(t,x)\leq W(t,x),\ 0\leq t<t+\delta \leq T,\ x\in \bar{D}$.\end{lemma}

\noindent{\it Proof}. For $u_1(\cdot)\in {\mathcal{U}}_{t, t+\delta}$ and
$u_2(\cdot)\in {\mathcal{U}}_{t+\delta, T}$, we define\ $u_1\oplus
u_2:=u_1\textbf{1}_{[t, t+\delta]}+u_2\textbf{1}_{(t+\delta, T]}, $\
which lies in ${\mathcal{U}}_{t, T}$.  Note that there exists a
sequence $\{u_i^1,\ i\geq 1\}\subset {\mathcal{U}}_{t, t+\delta}$\
such that
$$W_\delta(t, x)=\esssup_{u_1 \in {\mathcal{U}}_{t,
t+\delta}}I_\delta(t, x, u_1)=\sup_{i\geq 1}I_\delta(t, x, u_i^1),\
\ \mbox{{\it P}-a.s.}$$ For any $\varepsilon>0$, we define $\widetilde{\Gamma}_i:=\{W_\delta(t, x)\leq I_\delta(t, x,
u_i^1)+\varepsilon\}\in {\mathcal{F}}_{t},\ \ i\geq 1.$\  Then the
following mutually disjoint events $\Gamma_1:=\widetilde{\Gamma}_1,\ \
\Gamma_i:=\widetilde{\Gamma}_i\backslash(\cup^{i-1}_{l=1}\widetilde{\Gamma}_l)\in
{\mathcal{F}}_{t},\ i\geq 2,$\ form a $(\Omega,
{\mathcal{F}}_{t})$-partition. It is obvious that $u^\varepsilon_1:=\sum_{i\geq 1}\textbf{1}_{\Gamma_i}u_i^1\in
{\mathcal{U}}_{t, t+\delta}.$\ Moreover, from the uniqueness of the
solution of the forward-backward SDE, we have $I_\delta(t,
x, u^\varepsilon_1)=\sum_{i\geq 1}\textbf{1}_{\Gamma_i}I_\delta(t,
x, u_i^1),\ \ \mbox{{\it P}-a.s.}$\ Hence,
\begin{eqnarray}\label{3.17}
W_\delta(t,x)&\leq& \sum_{i\geq 1}\textbf{1}_{\Gamma_i}I_\delta(t,
x, u_i^1) +\varepsilon=I_\delta(t, x, u^\varepsilon_1)+\varepsilon\\
\nonumber &=& G^{t,x;u^\varepsilon_1}_{t,t+\delta} [W(t+\delta,
X^{t,x;u^\varepsilon_1}_{t+\delta})]+\varepsilon,\ \mbox{{\it
P}-a.s.}
\end{eqnarray}
 On the other hand,  from
the definition of $W(t+\delta,y)$\ we have, for any $y\in \bar{D},$
$$W(t+\delta,y)=\esssup_{u_2 \in {\mathcal{U}}_{t+\delta, T}}J(t+\delta, y; u_2),
\quad \mbox{{\it P}-a.s.}$$ Finally, since there exists a constant
$C\in {\mathbb{R}}$ such that for any $y, y' \in \bar{D},$\ $u_2\in {\mathcal{U}}_{t+\delta, T}$,
\be\label{3.18}
\begin{array}{llll}
{\rm(i)} & |W(t+\delta,y)-W(t+\delta,y')| \leq C(|y-y'|+|y-y'|^{\frac{1}{2}});  \\
{\rm(ii)} & |J(t+\delta, y, u_2)-J(t+\delta, y',
u_2)| \leq C(|y-y'|+|y-y'|^{\frac{1}{2}}),\ \ \mbox{{\it P}-a.s.,}\\
\end{array}
\ee (see Lemma \ref{l3.2}(i) and (\ref{3.6})(i)) we can prove by
approximating $X^{t,x;u_1^\varepsilon}_{t+\delta}$\ that
$$W(t+\delta, X^{t,x;u_1^\varepsilon}_{t+\delta} )\leq
\esssup_{u_2 \in {\mathcal{U}}_{t+\delta, T}}J(t+\delta,
X^{t,x;u_1^\varepsilon}_{t+\delta}; u_2),\ \mbox{{\it P}-a.s.}$$ To
estimate the right side of the above inequality we notice that there
exists some sequence $\{u_j^2,\ j\geq 1\}\subset
{\mathcal{U}}_{t+\delta, T}$\ such that
\begin{eqnarray*}
\esssup_{u_2 \in {\mathcal{U}}_{t+\delta,
T}}J(t+\delta,X^{t,x;u_1^\varepsilon}_{t+\delta}; u_2)=\sup_{j\geq
1}J(t+\delta,X^{t,x;u_1^\varepsilon}_{t+\delta}; u^2_j),\ \mbox{{\it
P}-a.s.}
\end{eqnarray*}
 Then, putting
$\widetilde{\Delta}_j:=\{\esssup_{u_2 \in {\mathcal{U}}_{t+\delta,
T}}J(t+\delta,X^{t,x;u_1^\varepsilon}_{t+\delta}; u_2)\leq
J(t+\delta,X^{t,x;u_1^\varepsilon}_{t+\delta};
u^2_j)+\varepsilon\}\in {\mathcal{F}}_{t+\delta},\ j\geq 1;$\ we
have with $\Delta_1:=\widetilde{\Delta}_1,\
\Delta_j:=\widetilde{\Delta}_j\backslash(\cup^{j-1}_{l=1}\widetilde{\Delta}_l)\in
{\mathcal{F}}_{t+\delta},\ j\geq 2,$\ an $(\Omega,
{\mathcal{F}}_{t+\delta})$-partition and
$u^\varepsilon_2:=\sum_{j\geq 1}\textbf{1}_{\Delta_j}u_j^2$\
 $\in {\mathcal{U}}_{t+\delta, T}.$ Therefore, from the uniqueness of the
solution of our reflected SDE and GBSDE, we have
\begin{eqnarray*}
& &J(t+\delta,X^{t,x;u_1^\varepsilon}_{t+\delta};
u_2^\varepsilon)=Y_{t+\delta}^{t+\delta,X^{t,x;u_1^\varepsilon}_{t+\delta};
u_2^\varepsilon}\ \ \ \ \hfill{\mbox{(see (\ref{3.8}))}}\\
&=&\sum_{j\geq
1}\textbf{1}_{\Delta_j}Y_{t+\delta}^{t+\delta,X^{t,x;u_1^\varepsilon}_{t+\delta};
 u_j^2}=\sum_{j\geq
1}\textbf{1}_{\Delta_j}J(t+\delta,X^{t,x;u_1^\varepsilon}_{t+\delta};
u_j^2),\ P{\mbox{-a.s}}.
\end{eqnarray*}
Thus,
\be\label{3.19}\begin{array}{lll}
& & W(t+\delta, X^{t,x;u_1^\varepsilon}_{t+\delta} )\leq \esssup_{u_2
\in {\mathcal{U}}_{t+\delta,
T}}J(t+\delta,X^{t,x;u_1^\varepsilon}_{t+\delta}; u_2)\\
&\leq& \sum_{j\geq
1}\textbf{1}_{\Delta_j}Y_{t+\delta}^{t,x;u_1^\varepsilon\oplus
u_j^2}+\varepsilon =Y_{t+\delta}^{t,x;u_1^\varepsilon\oplus u^\varepsilon_2}+\varepsilon
=Y_{t+\delta}^{t,x;u^\varepsilon}+\varepsilon,\
 \mbox{{\it P}-a.s.,}
\end{array}
\ee
where $u^\varepsilon:= u_1^\varepsilon\oplus u^\varepsilon_2\in
{\mathcal{U}}_{t, T}.$\ From (\ref{3.17}) and (\ref{3.19}) and
Lemmas 2.4 and 2.3, we get
 \be\label{3.20}
\begin{array}{lll}
  W_\delta(t,x)&\leq& G^{t,x;u^\varepsilon_1}_{t,t+\delta}
[Y_{t+\delta}^{t,x;u^\varepsilon}+\varepsilon]+\varepsilon \leq G^{t,x;u^\varepsilon_1}_{t,t+\delta}
[Y_{t+\delta}^{t,x;u^\varepsilon}]+
(C+1)\varepsilon\\
& =& G^{t,x;u^\varepsilon}_{t,t+\delta}
[Y_{t+\delta}^{t,x;u^\varepsilon}]+
(C+1)\varepsilon = Y_{t}^{t,x;u^\varepsilon}+
(C+1)\varepsilon\\
&\leq &  \esssup_{u \in {\mathcal{U}}_{t,
T}}Y_{t}^{t,x;u}+ (C+1)\varepsilon,\ \mbox{{\it P}-a.s.}
\end{array}
\ee That is,
 \be\label{3.21} W_\delta(t,x)\leq  W(t, x)+
(C+1)\varepsilon.\ee Finally, letting $\varepsilon\downarrow0,\
\mbox{we get}\ W_\delta(t,x)\leq W(t, x).$\endproof

\begin{lemma}\label{l3.5}$ W(t, x)\leq W_\delta(t,x),\ 0\leq t<t+\delta \leq T,\ x\in \bar{D}.$\end{lemma}
\unskip

\noindent{\it Proof}. Since $ W_\delta(t,x)=\esssup_{u_1 \in {\mathcal{U}}_{t,
t+\delta}}I_\delta(t, x, u_1), $\ \
we have \be
\begin{array}{lll}\label{3.22}
W_\delta(t,x)\geq I_\delta(t, x, u_1)=G^{t,x;u_1}_{t,t+\delta} [W(t+\delta, X^{t,x;u_1}_{t+\delta})],\
\end{array}\ee
$\mbox{{\it P}-a.s., for all}\ \ u_1\in {\mathcal{U}}_{t, t+\delta}$. Moreover, from the definition of $W(t+\delta,y),\ y\in \bar{D},$ we get
 \be\label{3.23}
W(t+\delta,y)= \esssup_{u_2 \in
{\mathcal{U}}_{t+\delta,T}}J(t+\delta, y; u_2),\ \mbox{{\it
P}-a.s.}\ee Let $\{O_i\}_{i\geq1}\subset
{\mathcal{B}}({\mathbb{R}}^d)$\ be a decomposition of
$\bar{D}$\ such that
$\sum\nolimits_{i\geq1}O_i=\bar{D}\ \mbox{and}\
\mbox{diam}(O_i)\break \leq \varepsilon,\ i\geq 1.$\ Let $y_i$\ be
an arbitrarily given element of $O_i,\ i\geq1.$\ We define $[X^{t,x;u_1}_{t+\delta}]:=\sum\nolimits_{i\geq1}y_i\textbf{1}_{\{X^{t,x;u_1}_{t+\delta}\in
\, O_i\}}.$\ Then we have \be\label{3.24}
|X^{t,x;u_1}_{t+\delta}-[X^{t,x;u_1}_{t+\delta}]|\leq \varepsilon,\
\mbox{everywhere on}\ \Omega, \ \mbox{for all}\ u_1\in
{\mathcal{U}}_{t, t+\delta}.\ee

Let $u\in {\mathcal{U}}_{t, T}$ be arbitrarily given and
decomposed into $u_1=u|_{[t, t+\delta]}\in {\mathcal{U}}_{t,
t+\delta}$\ and $u_2=u|_{(t+\delta, T]}\in {\mathcal{U}}_{t+\delta,
T}.$\ Then, from (\ref{3.22}), (\ref{3.18})(i), (\ref{3.24}), and
Lemmas 2.4 and 2.3, we have{\allowdisplaybreaks
\be\begin{array}{lll}\label{3.25}
W_\delta(t,x)&\geq&  G^{t,x;u_1}_{t,t+\delta} [W(t+\delta,
X^{t,x;u_1}_{t+\delta})] \geq G^{t,x;u_1}_{t,t+\delta}[W(t+\delta, [X^{t,x;u_1}_{t+\delta}])
-C\varepsilon-C\varepsilon^{\frac{1}{2}}]-\varepsilon\\
 &\geq& G^{t,x;u_1}_{t,t+\delta}[W(t+\delta,
[X^{t,x;u_1}_{t+\delta}])]-
C\varepsilon-C'\varepsilon^{\frac{1}{2}}\\  &=&
G^{t,x;u_1}_{t,t+\delta}[\sum\limits_{i\geq1}\textbf{1}_{\{X^{t,x;u_1}_{t+\delta}\in
O_i\}}W(t+\delta,y_i)]-C\varepsilon-C'\varepsilon^{\frac{1}{2}},\ \ \mbox{{\it P}-a.s.}
\end{array}
\ee
\noindent Furthermore, from} (\ref{3.23}), (\ref{3.18})(ii),
(\ref{3.24}), and Lemmas 2.4 and 2.3,
\begin{eqnarray}\label{3.26}
W_\delta(t,x)
&\geq&G^{t,x;u_1}_{t,t+\delta}[\sum\limits_{i\geq1}\textbf{1}_{\{X^{t,x;u_1}_{t+\delta}\in
O_i\}}J(t+\delta, y_i;
u_2)]-C\varepsilon-C'\varepsilon^{\frac{1}{2}}\\ \nonumber
&=&G^{t,x;u_1}_{t,t+\delta}[J(t+\delta,
[X^{t,x;u_1}_{t+\delta}];
u_2)]-C\varepsilon-C'\varepsilon^{\frac{1}{2}}\\ \nonumber
&\geq &G^{t,x;u_1}_{t,t+\delta}[J(t+\delta,X^{t,x;u_1}_{t+\delta};
u_2)-C''\varepsilon-C''\varepsilon^{\frac{1}{2}}]- C\varepsilon-C'\varepsilon^{\frac{1}{2}}\\
\nonumber &\geq
&G^{t,x;u_1}_{t,t+\delta}[J(t+\delta,X^{t,x;u_1}_{t+\delta}; u_2)]-
C\varepsilon-C'\varepsilon^{\frac{1}{2}}\\ \nonumber &=&
G^{t,x;u}_{t,t+\delta}[Y_{t+\delta}^{t,
x, u}]- C\varepsilon-C'\varepsilon^{\frac{1}{2}}\\
&=& Y_{t}^{t, x; u}- C\varepsilon-C'\varepsilon^{\frac{1}{2}},\
\mbox{{\it P}-a.s., for any}\ u\in {\mathcal{U}}_{t, T},\nonumber
\end{eqnarray}
where the constants $C, C', C''$\ may vary from lines to lines.
Consequently, \be\label{3.27}
\begin{array}{llll}
W_\delta(t,x)\geq\esssup_{u \in {\mathcal{U}}_{t,
T}}J(t, x; u)- C\varepsilon-C'\varepsilon^{\frac{1}{2}}=W(t,x)- C\varepsilon-C'\varepsilon^{\frac{1}{2}},\ \mbox{{\it
P}-a.s.}
\end{array}
\ee Finally, letting $\varepsilon\downarrow0$\ we get
$W_\delta(t,x)\geq W(t,x).$\ The proof is complete.\qquad\endproof

\begin{remark}\label{r3.4} {\rm (i)}  For any $u\in {\cal{U}}_{t, t+\delta},$
\be\label{3.29} W(t,x)(=W_\delta(t, x))\geq G^{t,x; u}_{t,t+\delta}
      [W(t+\delta, X^{t,x;u}_{t+\delta})],\quad \mbox{{\it P}-a.s.}
\ee

{\rm (ii)} From the inequality (\ref{3.17}), for all $(t, x)\in
[0,T]\times {\mathbb{R}}^n,$\ $\delta \in (0, T-t]$\ and
$\varepsilon>0$, the following holds:
 there exists some $u^{\varepsilon}(\cdot) \in {\cal{U}}_{t,
t+\delta}$\ such that
 \be\label{3.28} W(t,x)(=W_\delta(t, x))\leq G^{t,x;
u^{\varepsilon}}_{t,t+\delta}
      [W(t+\delta, X^{t,x; u^{\varepsilon}}_{t+\delta})]+C\varepsilon,\ \mbox{{\it P}-a.s.}
\ee

{\rm (iii)} Recall that the value function $W$ is deterministic. Then,
with $\delta=T-t$ and taking the expectation on both sides of (\ref{3.29})
 and (\ref{3.28}) we can get that $$ W(t,x)= \sup_{u
\in {\mathcal{U}}_{t,T}}E[J(t,x; u)]. $$
\end{remark}

\end{document}